# Mathematical Properties of Polynomial Dimensional Decomposition [*]

Sharif Rahman[†]

**Abstract.** Many high-dimensional uncertainty quantification problems are solved by polynomial dimensional decomposition (PDD), which represents Fourier-like series expansion in terms of random orthonormal polynomials with increasing dimensions. This study constructs dimension-wise and orthogonal splitting of polynomial spaces, proves completeness of polynomial orthogonal basis for prescribed assumptions, and demonstrates mean-square convergence to the correct limit – all associated with PDD. A second-moment error analysis reveals that PDD cannot commit larger error than polynomial chaos expansion (PCE) for the appropriately chosen truncation parameters. From the comparison of computational efforts, required to estimate with the same precision the variance of an output function involving exponentially attenuating expansion coefficients, the PDD approximation can be markedly more efficient than the PCE approximation.

**1. Introduction.** Polynomial dimensional decomposition (PDD) is a hierarchical, infinite series expansion of a square-integrable random variable involving measure-consistent orthogonal polynomials in independent random variables. Introduced by the author [20, 21] as a polynomial variant of the well-known analysis-of-variance (ANOVA) dimensional decomposition (ADD), PDD deflates the curse of dimensionality to some extent by developing an input-output behavior of complex systems with low effective dimensions [4], wherein the effects of degrees of interactions among input variables weaken rapidly or vanish altogether. Approximations stemming from truncated PDD are commonly used for solving uncertainty quantification problems in engineering and applied sciences, including multiscale fracture mechanics [6], random eigenvalue problems [24], computational fluid dynamics [27], and stochastic design optimization [25], to name a few. However, all existing works on PDD are focused on practical applications with almost no mathematical analysis of PDD. Indeed, a number of mathematical issues concerning necessary and sufficient conditions for the completeness of PDD basis functions; convergence, exactness, and optimal analyses of PDD; and approximation quality of the truncated PDD have yet to be studied or resolved. This paper fills the gap by establishing fundamental mathematical properties to empower PDD with a solid foundation, so that PDD can be as credible as its close cousin, polynomial chaos expansion (PCE) [5, 10, 28, 29], providing an alternative, if not a better, choice for uncertainty quantification in computational science and engineering.

The principal objective of this work is to examine important mathematical properties of PDD, not studied heretofore, for arbitrary but independent probability measures of input random variables. The paper is organized as follows. Section 2 defines or discusses mathematical notations and preliminaries. Two sets of assumptions on the input probability measures required by PDD are explained. A brief exposition of univariate and multivariate orthogonal polynomials consistent with a general but product-type probability measure, including their second moment properties, is given in Section 3. The section also describes relevant polynomial spaces and construction of their dimension-wise orthogonal decompositions. The orthogonal basis and completeness of multivariate orthogonal polynomials have also been proved. Section 4 briefly explains ADD, followed by presentations of PDD for a square-integrable random variable. The convergence and exactness of PDD are explained. In the same section, a truncated PDD and its approximation quality are discussed. The formulae for the mean and variance of a truncated PDD are also derived. The section ends with an explanation on how and when the PDD can be extended for infinitely many input variables. Section 5 briefly describes degree-wise orthogonal decompositions of polynomial spaces, leading to PCE. In Section 6, a second-moment error analysis of PDD is conducted, followed by a comparison with that of PCE. Finally, conclusions are drawn in Section 7.

---

[*]This work was supported by the U.S. National Science Foundation under Grant Number CMMI-1462385.
[†]College of Engineering and Applied Mathematics & Computational Sciences, The University of Iowa, Iowa City, IA 52242 (sharif-rahman@uiowa.edu). Questions, comments, or corrections to this document may be directed to that email address.

**2. Input Random Variables.** Let $\mathbb{N} := \{1, 2, \ldots\}$, $\mathbb{N}_0 := \mathbb{N} \cup \{0\}$, and $\mathbb{R} := (-\infty, +\infty)$ represent the sets of positive integer (natural), non-negative integer, and real numbers, respectively. Denote by $\mathbb{A}^{\{i\}}$, $i = 1, \ldots, N$, an $i$th bounded or unbounded subdomain of $\mathbb{R}$, so that $\mathbb{A}^N := \times_{i=1}^N \mathbb{A}^{\{i\}} \subseteq \mathbb{R}^N$.

Let $(\Omega, \mathcal{F}, \mathbb{P})$ be a complete probability space, where $\Omega$ is a sample space representing an abstract set of elementary events, $\mathcal{F}$ is a $\sigma$-algebra on $\Omega$, and $\mathbb{P} : \mathcal{F} \to [0, 1]$ is a probability measure. With $\mathcal{B}^N := \mathcal{B}(\mathbb{A}^N)$ representing the Borel $\sigma$-algebra on $\mathbb{A}^N \subseteq \mathbb{R}^N$, consider an $\mathbb{A}^N$-valued input random vector $\mathbf{X} := (X_1, \ldots, X_N)^T : (\Omega, \mathcal{F}) \to (\mathbb{A}^N, \mathcal{B}^N)$, describing the statistical uncertainties in all system parameters of a stochastic problem. The input random variables are also referred to as basic random variables [10]. The non-zero, finite integer $N$ represents the number of input random variables and is often referred to as the dimension of the stochastic problem.

Denote by $F_\mathbf{X}(\mathbf{x}) := \mathbb{P}(\cap_{i=1}^N \{X_i \leq x_i\})$ the joint distribution function of $\mathbf{X}$, admitting the joint probability density function $f_\mathbf{X}(\mathbf{x}) := \partial^N F_\mathbf{X}(\mathbf{x})/\partial x_1 \cdots \partial x_N$. Given the abstract probability space $(\Omega, \mathcal{F}, \mathbb{P})$ of $\mathbf{X}$, the image probability space is $(\mathbb{A}^N, \mathcal{B}^N, f_\mathbf{X} d\mathbf{x})$, where $\mathbb{A}^N$ can be viewed as the image of $\Omega$ from the mapping $\mathbf{X} : \Omega \to \mathbb{A}^N$, and is also the support of $f_\mathbf{X}(\mathbf{x})$. Similarly, each component random variable $X_i$ is defined on the abstract marginal probability space $(\Omega^{\{i\}}, \mathcal{F}^{\{i\}}, \mathbb{P}^{\{i\}})$ comprising sample space $\Omega^{\{i\}}$, $\sigma$-algebra $\mathcal{F}^{\{i\}}$, and probability measure $\mathbb{P}^{\{i\}}$. Then, the corresponding image probability space is $(\mathbb{A}^{\{i\}}, \mathcal{B}^{\{i\}}, f_{X_i} dx_i)$, where $\mathbb{A}^{\{i\}} \subseteq \mathbb{R}$ is the image sample space of $X_i$, $\mathcal{B}^{\{i\}}$ is the Borel $\sigma$-algebra on $\mathbb{A}^{\{i\}}$, and $f_{X_i}(x_i)$ is the marginal probability density function of $X_i$. Relevant statements and objects in the abstract probability space have obvious counterparts in the associated image probability space. Both probability spaces will be used in this paper.

Two sets of assumptions used by PDD are as follows.

**Assumption 2.1.** *The input random vector* $\mathbf{X} := (X_1, \ldots, X_N)^T : (\Omega, \mathcal{F}) \to (\mathbb{A}^N, \mathcal{B}^N)$ *satisfies all of the following conditions:*

(1) *Each input random variable* $X_i : (\Omega^{\{i\}}, \mathcal{F}^{\{i\}}) \to (\mathbb{A}^{\{i\}}, \mathcal{B}^{\{i\}})$ *has absolutely continuous marginal distribution function* $F_{X_i}(x_i) := \mathbb{P}(X_i \leq x_i)$ *and continuous marginal density function* $f_{X_i}(x_i) := \partial F_{X_i}(x_i)/\partial x_i$ *with a bounded or unbounded support* $\mathbb{A}^{\{i\}} \subseteq \mathbb{R}$.

(2) *All component random variables* $X_i$, $i = 1, \ldots, N$, *are statistically independent, but not necessarily identical. In consequence,* $\mathbf{X}$ *is endowed with a product-type probability density function, that is,* $f_\mathbf{X}(\mathbf{x}) = \prod_{i=1}^N f_{X_i}(x_i)$ *with a bounded or unbounded support* $\mathbb{A}^N \subseteq \mathbb{R}^N$.

(3) *Each input random variable* $X_i$ *possesses finite moments of all orders, that is, for all*

$i = 1, \ldots, N$ and $l \in \mathbb{N}_0$,

$$(2.1) \quad \mu_{i,l} := \mathbb{E}\left[X_i^l\right] := \int_\Omega X_i^l(\omega) d\mathbb{P}(\omega) = \int_{\mathbb{A}^N} x_i^l f_{\mathbf{X}}(\mathbf{x}) d\mathbf{x} = \int_{\mathbb{A}^{\{i\}}} x_i^l f_{X_i}(x_i) dx_i < \infty,$$

where $\mathbb{E}$ is the expectation operator with respect to the probability measure $\mathbb{P}$ or $f_{\mathbf{X}}(\mathbf{x})d\mathbf{x}$.

**Assumption 2.2.** *The moments and marginal density function of each input random variable $X_i$, where $i = 1, \ldots, N$, satisfy at least one of the following conditions [10]:*
(1) *The density function $f_{X_i}(x_i)$ has a compact support, that is, there exists a compact interval $[a_i, b_i]$, $a_i, b_i \in \mathbb{R}$, such that $\mathbb{P}(a_i \leq X_i \leq b_i) = 1$.*
(2) *For the moment sequence $\{\mu_{i,l}\}_{l \in \mathbb{N}_0}$ for $X_i$, there holds*

$$(2.2) \quad \liminf_{l \to \infty} \frac{(\mu_{i,2l})^{1/2l}}{2l} < \infty.$$

(3) *For the moment sequence $\{\mu_{i,l}\}_{l \in \mathbb{N}_0}$ for $X_i$, there holds*

$$(2.3) \quad \sum_{l \in \mathbb{N}_0} \frac{1}{(\mu_{i,2l})^{1/2l}} = \infty.$$

(4) *The random variable $X_i$ is exponentially integrable, that is, there exists a real number $a > 0$ such that*

$$(2.4) \quad \int_{\mathbb{A}^{\{i\}}} \exp(a|x_i|) f_{X_i}(x_i) dx_i < \infty.$$

(5) *If the density function $f_{X_i}(x_i)$ is symmetric, differentiable, and strictly positive, then there exists a real number $a > 0$ such that*

$$(2.5) \quad \int_{\mathbb{A}^{\{i\}}} -\frac{\ln f_{X_i}(x_i)}{1 + x_i^2} dx_i = \infty \text{ and } \frac{-x_i df_{X_i}(x_i)/dx_i}{f_{X_i}(x_i)} \to \infty \text{ as } (x_i \to \infty, x_i \geq a).$$

Assumption 2.1 assures the existence of an infinite sequence of orthogonal polynomials consistent with the input probability measure. Assumption 2.2, in addition to Assumption 2.1, guarantees the input probability measure to be determinate, resulting in a complete orthogonal polynomial basis of a function space of interest. Both assumptions impose only mild restrictions on the probability measure. Examples of input random variables satisfying Assumptions 2.1 and 2.2 are Gaussian, uniform, exponential, beta, and gamma variables, which are commonly used in uncertainty quantification. These assumptions, to be explained in the next section, are vitally important for the determinacy of the probability measure and the completeness of the orthogonal polynomial basis. Therefore, for both PDD and PCE, which entail orthogonal polynomial expansions, Assumptions 2.1 and 2.2 are necessary. Unfortunately, they are not always clearly specified in the PDD or PCE literature. A prototypical example where Assumption 2.1 is satisfied, but Assumption 2.2 is not, is the case of a lognormal random variable. As noted by Ernst *et al.* [10], the violation of Assumption 2.2 leads to indeterminacy of the input probability measure and thereby fails to form a complete orthogonal polynomial basis. Finally, Assumptions 2.1 and 2.2 can be modified to account for random variables with discrete or mixed distributions [11] or dependent random variables [23]. The discrete or mixed distributions and dependent variables are not considered in this paper.

### 3. Measure-Consistent Orthogonal Polynomials and Polynomial Spaces.

**3.1. Univariate orthogonal polynomials.** Consider an $i$th random variable $X_i$ defined on the abstract probability space $(\Omega^{\{i\}}, \mathcal{F}^{\{i\}}, \mathbb{P}^{\{i\}})$ with its image $(\mathbb{A}^{\{i\}}, \mathcal{B}^{\{i\}}, f_{X_i} dx_i)$. Let $\Pi^{\{i\}} := \mathbb{R}[x_i]$ be the space of real polynomials in $x_i$. For any polynomial pair $P_{\{i\}}, Q_{\{i\}} \in \Pi^{\{i\}}$, define an inner product

$$(3.1) \quad (P_{\{i\}}, Q_{\{i\}})_{f_{X_i} dx_i} := \int_{\mathbb{A}^{\{i\}}} P_{\{i\}}(x_i) Q_{\{i\}}(x_i) f_{X_i}(x_i) dx_i =: \mathbb{E}\left[P_{\{i\}}(X_i) Q_{\{i\}}(X_i)\right]$$

with respect to the probability measure $f_{X_i}(x_i) dx_i$ and the induced norm

$$\|P_{\{i\}}\|_{f_{X_i} dx_i} := \sqrt{(P_{\{i\}}, P_{\{i\}})_{f_{X_i} dx_i}} = \left(\int_{\mathbb{A}_i} P_{\{i\}}^2(x_i) f_{X_i}(x_i) dx_i\right)^{1/2} =: \sqrt{\mathbb{E}\left[P_{\{i\}}^2(X_i)\right]}.$$

Under Assumption 2.1, moments of $X_i$ of all orders exist and are finite, including *zero-order moments* $\mu_{i,0} := \int_{\mathbb{A}^{\{i\}}} f_{X_i}(x_i) dx_i = 1$, $i = 1, \ldots, N$, that are always positive. Clearly, $\|P_{\{i\}}\| > 0$ for all non-zero $P_{\{i\}} \in \Pi^{\{i\}}$. Then, according to Gautschi [12], the inner product in (3.1) is positive-definite on $\Pi^{\{i\}}$. Therefore, there exists an infinite set of univariate orthogonal polynomials, say, $\{P_{\{i\},j_i}(x_i) : j_i \in \mathbb{N}_0\}$, $P_{\{i\},j_i} \neq 0$, which is consistent with the probability measure $f_{X_i}(x_i) dx_i$, satisfying

$$(3.2) \quad \left(P_{\{i\},j_i}, P_{\{i\},k_i}\right)_{f_{X_i} dx_i} = \begin{cases} \mathbb{E}[P_{\{i\},j_i}^2(X_i)], & j_i = k_i, \\ 0, & j_i \neq k_i, \end{cases}$$

for $k_i \in \mathbb{N}_0$, where $0 < \mathbb{E}[P_{\{i\},j_i}^2(X_i)] < \infty$. Here, in the notation for the polynomial $P_{\{i\},j_i}(x_i)$, the first and second indices refer to the $i$th variable and degree $j_i$, respectively. Prominent examples of classical univariate orthogonal polynomials comprise Hermite, Laguerre, and Jacobi polynomials, which are consistent with the measures defined by Gaussian, gamma, and beta densities on the whole real line, semi-infinite interval, and bounded interval, respectively. Many orthogonal polynomials, including the three classical polynomials mentioned, can be expressed in a unified way by invoking hypergeometric series, incorporated in a tree structure of the Askey scheme [1]. For even more general measures, established numerical techniques, such as Gram-Schmidt [13] and Stieltjes' procedure [26], can be used to generate any measure-consistent orthogonal polynomials.

**3.2. Multivariate orthogonal polynomials.** For $N \in \mathbb{N}$, denote by $\{1, \ldots, N\}$ an index set, so that $u \subseteq \{1, \ldots, N\}$ is a subset, including the empty set $\emptyset$, with cardinality $0 \leq |u| \leq N$. When $\emptyset \neq u \subseteq \{1, \ldots, N\}$, a $|u|$-dimensional multi-index is denoted by $\mathbf{j}_u := (j_{i_1}, \ldots, j_{i_{|u|}}) \in \mathbb{N}_0^{|u|}$ with degree $|\mathbf{j}_u| := j_{i_1} + \cdots + j_{i_{|u|}}$, where $j_{i_p} \in \mathbb{N}_0$, $p = 1, \ldots, |u|$, represents the $p$th component of $\mathbf{j}_u$.[1]

For $\emptyset \neq u \subseteq \{1, \ldots, N\}$, let $\mathbf{X}_u := (X_{i_1}, \ldots, X_{i_{|u|}})^T$, a subvector of $\mathbf{X}$, be defined on the abstract probability space $(\Omega^u, \mathcal{F}^u, \mathbb{P}^u)$, where $\Omega^u$ is the sample space of $\mathbf{X}_u$, $\mathcal{F}^u$ is a $\sigma$-algebra on $\Omega^u$, and $\mathbb{P}^u$ is a probability measure. The corresponding image probability space is $(\mathbb{A}^u, \mathcal{B}^u, f_{\mathbf{X}_u} d\mathbf{x}_u)$, where $\mathbb{A}^u := \times_{i \in u} \mathbb{A}^{\{i\}} \subseteq \mathbb{R}^{|u|}$ is the image sample space of $\mathbf{X}_u$, $\mathcal{B}^u$ is the Borel $\sigma$-algebra on $\mathbb{A}^u$, and $f_{\mathbf{X}_u}(\mathbf{x}_u)$ is the marginal probability density function of $\mathbf{X}_u$ supported on $\mathbb{A}^u$. Under Assumption 2.1, $f_{\mathbf{X}_u}(\mathbf{x}_u) = \prod_{i \in u} f_{X_i}(x_i)$. Denote by

$$\Pi^u := \mathbb{R}[\mathbf{x}_u] = \mathbb{R}[x_{i_1}, \ldots, x_{i_{|u|}}]$$

the space of all real polynomials in $\mathbf{x}_u$. Then, given the inner product

$$(P_u, Q_u)_{f_{\mathbf{X}_u} d\mathbf{x}_u} := \int_{\mathbb{A}^u} P_u(\mathbf{x}_u) Q_u(\mathbf{x}_u) f_{\mathbf{X}_u}(\mathbf{x}_u) d\mathbf{x}_u =: \mathbb{E}\left[P_u(\mathbf{X}_u) Q_u(\mathbf{X}_u)\right],$$

two polynomials $P_u \in \Pi^u$ and $Q_u \in \Pi^u$ in $\mathbf{x}_u$ are called orthogonal to each other if $(P_u, Q_u)_{f_{\mathbf{X}_u} d\mathbf{x}_u} = 0$ [8]. Moreover, a polynomial $P_u \in \Pi^u$ is said to be an orthogonal polynomial with respect to $f_{\mathbf{X}_u} d\mathbf{x}_u$ if it is orthogonal to all polynomials of lower degree, that is, if [8]

$$(3.3) \qquad (P_u, Q_u)_{f_{\mathbf{X}_u} d\mathbf{x}_u} = 0 \; \forall Q_u \in \Pi^u \text{ with } \deg Q_u < \deg P_u.$$

---

[1] The same symbol $|\cdot|$ is used for designating both the cardinality of a set and the degree of a multi-index in this paper.

Let $\{P_{u,\mathbf{j}_u}(\mathbf{x}_u) : \mathbf{j}_u \in \mathbb{N}_0^{|u|}\}$, $\emptyset \neq u \subseteq \{1, \ldots, N\}$, represent an infinite set of multivariate orthogonal polynomials, which is consistent with the probability measure $f_{\mathbf{X}_u}(\mathbf{x}_u) d\mathbf{x}_u$, satisfying

$$(3.4) \qquad (P_{u,\mathbf{j}_u}, P_{u,\mathbf{k}_u})_{f_{\mathbf{X}_u} d\mathbf{x}_u} =: \mathbb{E}\left[P_{u,\mathbf{j}_u}(\mathbf{X}_u) P_{u,\mathbf{k}_u}(\mathbf{X}_u)\right] = 0 \; \forall \mathbf{j}_u \neq \mathbf{k}_u, \; \mathbf{k}_u \in \mathbb{N}_0^{|u|}.$$

Clearly, each $P_{u,\mathbf{j}_u} \in \Pi^u$ is a multivariate orthogonal polynomial satisfying (3.3). Due to the product-type probability measure of $\mathbf{X}_u$, a consequence of statistical independence from Assumption 2.1, such multivariate polynomials exist and are easily constructed by tensorizing univariate orthogonal polynomials.

Proposition 3.1. *Let* $\mathbf{X} := (X_1, \ldots, X_N)^T : (\Omega, \mathcal{F}) \to (\mathbb{A}^N, \mathcal{B}^N)$ *be a vector of* $N \in \mathbb{N}$ *input random variables fulfilling Assumption 2.1. Suppose that the sets of univariate orthogonal polynomials for all marginal measures are obtained as* $\{P_{\{i\},j_i}(x_i) : j_i \in \mathbb{N}_0\}$, $i = 1, \ldots, N$. *Then, for* $\emptyset \neq u \subseteq \{1, \ldots, N\}$, *the set of multivariate orthogonal polynomials in* $\mathbf{x}_u$ *consistent with the probability measure* $f_{\mathbf{X}_u}(\mathbf{x}_u) d\mathbf{x}_u$ *is*

$$(3.5) \qquad \left\{P_{u,\mathbf{j}_u}(\mathbf{x}_u) : \mathbf{j}_u \in \mathbb{N}_0^{|u|}\right\} = \bigotimes_{i \in u} \left\{P_{\{i\},j_i}(x_i) : j_i \in \mathbb{N}_0\right\},$$

*where the symbol* $\bigotimes$ *denotes tensor product. In terms of an element, the multivariate orthogonal polynomial of degree* $|\mathbf{j}_u| = j_{i_1} + \cdots + j_{i_{|u|}}$ *is*

$$(3.6) \qquad P_{u,\mathbf{j}_u}(\mathbf{x}_u) = \prod_{i \in u} P_{\{i\},j_i}(x_i).$$

*Proof.* Consider two distinct polynomials $P_{u,\mathbf{j}_u}(\mathbf{x}_u)$ and $P_{u,\mathbf{k}_u}(\mathbf{x}_u)$ from the set $\{P_{u,\mathbf{j}_u}(\mathbf{x}_u) : \mathbf{j}_u \in \mathbb{N}_0^{|u|}\}$ satisfying (3.5). Since $\mathbf{j}_u \neq \mathbf{k}_u$, $\mathbf{j}_u$ and $\mathbf{k}_u$ must differ in at least one component. Without loss of generality, suppose that $j_{i_1} \neq k_{i_1}$. Then, by Fubini's theorem, with statistical independence of random variables in mind,

$$(3.7) \quad \begin{aligned} (P_{u,\mathbf{j}_u} P_{u,\mathbf{k}_u})_{f_{\mathbf{X}_u} d\mathbf{x}_u} &= \int_{\mathbb{A}^u} P_{u,\mathbf{j}_u}(\mathbf{x}_u) P_{u,\mathbf{k}_u}(\mathbf{x}_u) f_{\mathbf{X}_u}(\mathbf{x}_u) d\mathbf{x}_u \\ &= \int_{\times_{p=2}^{|u|} \mathbb{A}^{\{i_p\}}} \prod_{p=2}^{|u|} \left[P_{\{i_p\},j_{i_p}}(x_{i_p}) P_{\{i_p\},k_{i_p}}(x_{i_p}) f_{X_{i_p}}(x_{i_p}) dx_{i_p}\right] \times \\ & \qquad \int_{\mathbb{A}^{\{i_1\}}} P_{\{i_1\},j_{i_1}}(x_{i_1}) P_{\{i_1\},k_{i_1}}(x_{i_1}) f_{X_{i_1}}(x_{i_1}) dx_{i_1} \\ &= 0, \end{aligned}$$

where the equality to *zero* in the last line results from the recognition that the inner integral vanishes by setting $i = i_1$ in (3.2).

In addition, for $\mathbf{j}_u \in \mathbb{N}_0^{|u|}$,

$$(3.8) \qquad (P_{u,\mathbf{j}_u} P_{u,\mathbf{j}_u})_{f_{\mathbf{X}_u} d\mathbf{x}_u} =: \mathbb{E}\left[P_{u,\mathbf{j}_u}^2(\mathbf{X}_u)\right] = \prod_{i \in u} \mathbb{E}\left[P_{\{i\},j_i}^2(X_i)\right] > 0$$

and is finite by virtue of the existence of the set of univariate orthogonal polynomials $\{P_{\{i\},j_i}(x_i) : j_i \in \mathbb{N}_0\}$ for $i = 1, \ldots, N$. Therefore, $\{P_{u,\mathbf{j}_u}(\mathbf{x}_u) : \mathbf{j}_u \in \mathbb{N}_0^{|u|}\}$ satisfying (3.5) is a set of multivariate orthogonal polynomials consistent with the probability measure $f_{\mathbf{X}_u}(\mathbf{x}_u) d\mathbf{x}_u$. ∎

Once the multivariate orthogonal polynomials are obtained, they can be scaled to generate multivariate orthonormal polynomials, as follows.

**Definition 3.2.** *A multivariate orthonormal polynomial* $\Psi_{u,\mathbf{j}_u}(\mathbf{x}_u)$, $\emptyset \neq u \subseteq \{1, \ldots, N\}$, $\mathbf{j}_u \in \mathbb{N}_0^{|u|}$, *of degree* $|\mathbf{j}_u| = j_{i_1} + \cdots + j_{i_{|u|}}$ *that is consistent with the probability measure* $f_{\mathbf{X}_u}(\mathbf{x}_u) d\mathbf{x}_u$ *is defined as*

$$(3.9) \qquad \Psi_{u,\mathbf{j}_u}(\mathbf{x}_u) := \frac{P_{u,\mathbf{j}_u}(\mathbf{x}_u)}{\sqrt{\mathbb{E}[P_{u,\mathbf{j}_u}^2(\mathbf{X}_u)]}} = \prod_{i \in u} \frac{P_{\{i\},j_i}(x_i)}{\sqrt{\mathbb{E}\left[P_{\{i\},j_i}^2(X_i)\right]}} =: \prod_{i \in u} \Psi_{\{i\},j_i}(x_i),$$

*where* $\Psi_{\{i\},j_i}(x_i) := P_{\{i\},j_i}(x_i) / \sqrt{\mathbb{E}[P_{\{i\},j_i}^2(X_i)]}$ *is a univariate orthonormal polynomial in* $x_i$ *of degree* $j_i$ *that is consistent with the probability measure* $f_{X_i}(x_i) dx_i$

### 3.3. Dimension-wise orthogonal decomposition of polynomial spaces.
An orthogonal decomposition of polynomial spaces entailing dimension-wise splitting leads to PDD. Here, to facilitate such splitting of the polynomial space $\Pi^u$ for any $\emptyset \neq u \subseteq \{1, \ldots, N\}$, limit the power $j_{i_p}$ of the $i_p$-th variable, where $i_p \in u \subseteq \{1, \ldots, N\}$, $p = 1, \ldots, |u|$, and $|u| > 0$, to take on only positive integer values. In consequence, $\mathbf{j}_u := (j_{i_1}, \ldots, j_{i_{|u|}}) \in \mathbb{N}^{|u|}$, the multi-index of $P_{u,\mathbf{j}_u}(\mathbf{x}_u)$, has degree $|\mathbf{j}_u| = j_{i_1} + \cdots + j_{i_{|u|}}$, varying from $|u|$ to $\infty$ as $j_{i_1} \neq \cdots j_{i_{|u|}} \neq 0$.

For $\mathbf{j}_u \in \mathbb{N}^{|u|}$ and $\mathbf{x}_u := (x_{i_1}, \ldots, x_{i_{|u|}})$, a monomial in the variables $x_{i_1}, \ldots, x_{i_{|u|}}$ is the product $\mathbf{x}_u^{\mathbf{j}_u} = x_{i_1}^{j_{i_1}} \ldots x_{i_{|u|}}^{j_{i_{|u|}}}$ and has a total degree $|\mathbf{j}_u|$. A linear combination of $\mathbf{x}_u^{\mathbf{j}_u}$, where $|\mathbf{j}_u| = l$, $|u| \leq l \leq \infty$, is a homogeneous polynomial in $\mathbf{x}_u$ of degree $l$. For $\emptyset \neq u \subseteq \{1, \ldots, N\}$, denote by

$$\mathcal{Q}_l^u := \mathrm{span}\{\mathbf{x}_u^{\mathbf{j}_u} : |\mathbf{j}_u| = l, \mathbf{j}_u \in \mathbb{N}^{|u|}\}, \ |u| \leq l < \infty,$$

the space of homogeneous polynomials in $\mathbf{x}_u$ of degree $l$ where the individual degree of each variable is non-zero and by

$$\Theta_m^u := \mathrm{span}\{\mathbf{x}_u^{\mathbf{j}_u} : |u| \leq |\mathbf{j}_u| \leq m, \mathbf{j}_u \in \mathbb{N}^{|u|}\}, \ |u| \leq m < \infty,$$

the space of polynomials in $\mathbf{x}_u$ of degree at least $|u|$ and at most $m$ where the individual degree of each variable is non-zero. The dimensions of the vector spaces $\mathcal{Q}_l^u$ and $\Theta_m^u$, respectively, are

$$(3.10) \qquad \dim \mathcal{Q}_l^u = \#\left\{\mathbf{j}_u \in \mathbb{N}^{|u|} : |\mathbf{j}_u| = l\right\} = \binom{l-1}{|u|-1}$$

and

$$\dim \Theta_m^u = \sum_{l=|u|}^m \dim \mathcal{Q}_l^u = \sum_{l=|u|}^m \binom{l-1}{|u|-1} = \binom{m}{|u|}. \tag{3.11}$$

Let $\mathcal{Z}_{|u|}^u := \Theta_{|u|}^u$. For each $|u|+1 \le l < \infty$, denote by $\mathcal{Z}_l^u \subset \Theta_l^u$ the space of orthogonal polynomials of degree exactly $l$ that are orthogonal to all polynomials in $\Theta_{l-1}^u$, that is,

$$\mathcal{Z}_l^u := \{P_u \in \Theta_l^u : (P_u, Q_u)_{f_{\mathbf{X}_u} d\mathbf{x}_u} = 0 \; \forall \, Q_u \in \Theta_{l-1}^u\}, \; |u|+1 \le l < \infty.$$

Then $\mathcal{Z}_l^u$, provided that the support of $f_{\mathbf{X}_u}(\mathbf{x}_u)$ has non-empty interior, is a vector space of dimension

$$M_{u,l} := \dim \mathcal{Z}_l^u = \dim \mathcal{Q}_l^u = \binom{l-1}{|u|-1}.$$

Many choices exist for the basis of $\mathcal{Z}_l^u$. Here, to be formally proved in Section 3.3.2, select $\{P_{u,\mathbf{j}_u}(\mathbf{x}_u) : |\mathbf{j}_u| = l, \mathbf{j}_u \in \mathbb{N}^{|u|}\} \subset \mathcal{Z}_l^u$ to be a basis of $\mathcal{Z}_l^u$, comprising $M_{u,l}$ number of basis functions. Each basis function $P_{u,\mathbf{j}_u}(\mathbf{x}_u)$ is a multivariate orthogonal polynomial of degree $|\mathbf{j}_u|$ as defined earlier. Clearly,

$$\mathcal{Z}_l^u = \mathrm{span}\{P_{u,\mathbf{j}_u}(\mathbf{x}_u) : |\mathbf{j}_u| = l, \mathbf{j}_u \in \mathbb{N}^{|u|}\}, \; |u| \le l < \infty.$$

According to Proposition 3.3, to be presented later, $P_{u,\mathbf{j}_u}(\mathbf{X}_u)$ is orthogonal to $P_{v,\mathbf{k}_v}(\mathbf{X}_v)$ whenever (1) $u \ne v$ and $\mathbf{j}_u, \mathbf{k}_v$ are arbitrary; or (2) $u = v$ and $\mathbf{j}_u \ne \mathbf{k}_v$. Therefore, any two distinct polynomial subspaces $\mathcal{Z}_l^u$ and $\mathcal{Z}_{l'}^v$, where $\emptyset \ne u \subseteq \{1,\ldots,N\}$, $\emptyset \ne v \subseteq \{1,\ldots,N\}$, $|u| \le l < \infty$, and $|v| \le l' < \infty$, are orthogonal whenever $u \ne v$ or $l \ne l'$. In consequence, there exist orthogonal decompositions of

$$\begin{aligned}
\Theta_m^u &= \bigoplus_{l=|u|}^m \mathcal{Z}_l^u = \bigoplus_{l=|u|}^m \mathrm{span}\{P_{u,\mathbf{j}_u}(\mathbf{x}_u) : |\mathbf{j}_u| = l, \mathbf{j}_u \in \mathbb{N}^{|u|}\} \\
&= \mathrm{span}\{P_{u,\mathbf{j}_u}(\mathbf{x}_u) : |u| \le |\mathbf{j}_u| \le m, \mathbf{j}_u \in \mathbb{N}^{|u|}\}
\end{aligned}$$

with the symbol $\oplus$ representing orthogonal sum and

$$\begin{aligned}
\Pi^u &= \mathbf{1} \oplus \bigoplus_{\emptyset \ne v \subseteq u} \bigoplus_{l=|v|}^\infty \mathcal{Z}_l^v = \mathbf{1} \oplus \bigoplus_{\emptyset \ne v \subseteq u} \bigoplus_{l=|v|}^\infty \mathrm{span}\{P_{v,\mathbf{j}_v}(\mathbf{x}_v) : |\mathbf{j}_v| = l, \mathbf{j}_v \in \mathbb{N}^{|v|}\} \\
&= \mathbf{1} \oplus \bigoplus_{\emptyset \ne v \subseteq u} \mathrm{span}\{P_{v,\mathbf{j}_v}(\mathbf{x}_v) : \mathbf{j}_v \in \mathbb{N}^{|v|}\},
\end{aligned} \tag{3.12}$$

where $\mathbf{1} := \mathrm{span}\{1\}$, the constant subspace, needs to be added because the subspace $\mathcal{Z}_l^v$ excludes constant functions.

Recall that $\Pi^N$ is the space of all real polynomials in $\mathbf{x}$. Then, setting $u = \{1,\ldots,N\}$ in (3.12) first and then swapping $v$ for $u$ yields yet another orthogonal decomposition of

$$\begin{aligned}
\Pi^N &= \mathbf{1} \oplus \bigoplus_{\emptyset \ne u \subseteq \{1,\ldots,N\}} \bigoplus_{l=|u|}^\infty \mathcal{Z}_l^u \\
&= \mathbf{1} \oplus \bigoplus_{\emptyset \ne u \subseteq \{1,\ldots,N\}} \bigoplus_{l=|u|}^\infty \mathrm{span}\{P_{u,\mathbf{j}_u}(\mathbf{x}_u) : |\mathbf{j}_u| = l, \mathbf{j}_u \in \mathbb{N}^{|u|}\} \\
&= \mathbf{1} \oplus \bigoplus_{\emptyset \ne u \subseteq \{1,\ldots,N\}} \mathrm{span}\{P_{u,\mathbf{j}_u}(\mathbf{x}_u) : \mathbf{j}_u \in \mathbb{N}^{|u|}\}.
\end{aligned} \tag{3.13}$$

Note that the last expression of (3.13) is equal to the span of

$$\left\{P_{\mathbf{j}}(\mathbf{x}) : \mathbf{j} \in \mathbb{N}_0^N\right\} := \bigotimes_{i=1}^{N} \left\{P_{\{i\},j_i}(x_i) : j_i \in \mathbb{N}_0\right\}, \quad (3.14)$$

representing an infinite set of orthogonal polynomials in $\mathbf{x}$.

Given the dimension-wise orthogonal splitting of $\Pi^N$, any square-integrable function of input random vector $\mathbf{X}$ can be expanded as a Fourier-like series of hierarchically ordered multivariate orthogonal or orthonormal polynomials in $\mathbf{X}_u$, $\emptyset \neq u \subseteq \{1, \ldots, N\}$. The expansion is referred to as PDD, to be formally presented and analyzed in Section 4.

**3.4. Statistical properties of random multivariate polynomials.** When the input random variables $X_1, \ldots, X_N$, instead of real variables $x_1, \ldots, x_N$, are inserted in the argument, the multivariate polynomials $P_{u,\mathbf{j}_u}(\mathbf{X}_u)$ and $\Psi_{u,\mathbf{j}_u}(\mathbf{X}_u)$, where $\emptyset \neq u \subseteq \{1, \ldots, N\}$ and $\mathbf{j}_u \in \mathbb{N}^{|u|}$, become functions of random input variables. Therefore, it is important to establish their second-moment properties, to be exploited in the remaining part of this section and Section 4.

**Proposition 3.3.** *Let $\mathbf{X} := (X_1, \ldots, X_N)$ be a vector of $N \in \mathbb{N}$ input random variables fulfilling Assumption 2.1. For $\emptyset \neq u, v \subseteq \{1, \ldots, N\}$, $\mathbf{j}_u \in \mathbb{N}^{|u|}$, and $\mathbf{k}_v \in \mathbb{N}^{|v|}$, the first- and second-order moments of multivariate orthogonal polynomials are*

$$\mathbb{E}\left[P_{u,\mathbf{j}_u}(\mathbf{X}_u)\right] = 0 \quad (3.15)$$

*and*

$$\mathbb{E}\left[P_{u,\mathbf{j}_u}(\mathbf{X}_u) P_{v,\mathbf{k}_v}(\mathbf{X}_v)\right] = \begin{cases} \prod_{i \in u} \mathbb{E}\left[P_{\{i\},j_i}^2(X_i)\right], & u = v, \ \mathbf{j}_u = \mathbf{k}_v, \\ 0, & \text{otherwise}, \end{cases} \quad (3.16)$$

*respectively.*

*Proof.* Using (3.6) and statistical independence of random variables, $\mathbb{E}[P_{u,\mathbf{j}_u}(\mathbf{X}_u)] = \prod_{i \in u} \mathbb{E}[P_{\{i\},j_i}(X_i)]$ for any $\mathbf{j}_u \in \mathbb{N}^{|u|}$. Since each component of $\mathbf{j}_u$ is non-zero, (3.2), with the constant function $P_{\{i\},0} \neq 0$ in mind, produces $\mathbb{E}[P_{\{i\},j_i}(X_i)] = 0$ for any $i \in u$, $j_i \in \mathbb{N}$, resulting in (3.15).

To obtain the non-trivial result of (3.16), set $u = v$ and $\mathbf{j}_u = \mathbf{k}_v$ and use (3.8) directly. The trivial result of (3.16) is obtained by considering two subcases. First, when $u = v$ and $\mathbf{j}_u \neq \mathbf{k}_v$, (3.7) yields the result already. Second, when $u \neq v$ and $\mathbf{j}_u, \mathbf{k}_v \in \mathbb{N}^{|v|}$ are arbitrary, then $u$ and $v$ differ by at least one element. Suppose that $i \in (u \cup v)$ is that element with the associated degree $j_i \in \mathbb{N}$. Using the statistical independence of random variables and the fact that $\mathbb{E}[P_{\{i\},j_i}(X_i)] = 0$, as already demonstrated, produces the desired result. ∎

**Corollary 3.4.** *For $\emptyset \neq u, v \subseteq \{1, \ldots, N\}$, $\mathbf{j}_u \in \mathbb{N}^{|u|}$, and $\mathbf{k}_v \in \mathbb{N}^{|v|}$, the first- and second-order moments of multivariate orthonormal polynomials are*

$$\mathbb{E}\left[\Psi_{u,\mathbf{j}_u}(\mathbf{X}_u)\right] = 0 \quad (3.17)$$

*and*

$$\mathbb{E}\left[\Psi_{u,\mathbf{j}_u}(\mathbf{X}_u) \Psi_{v,\mathbf{k}_v}(\mathbf{X}_v)\right] = \begin{cases} 1, & u = v, \ \mathbf{j}_u = \mathbf{k}_v, \\ 0, & \text{otherwise}, \end{cases} \quad (3.18)$$

*respectively.*

**3.5. Orthogonal basis and completeness.** An important question regarding multivariate orthogonal polynomials discussed in the preceding subsection is whether they constitute a complete basis in a function space of interest, such as a Hilbert space. Let $L^2(\mathbb{A}^N, \mathcal{B}^N, f_{\mathbf{X}}d\mathbf{x})$ represent a Hilbert space of square-integrable functions with respect to the probability measure $f_{\mathbf{X}}(\mathbf{x})d\mathbf{x}$ supported on $\mathbb{A}^N$. The following two propositions show that, indeed, measure-consistent orthogonal polynomials span various spaces of interest.

*Proposition 3.5.* Let $\mathbf{X} := (X_1, \ldots, X_N)^T : (\Omega, \mathcal{F}) \to (\mathbb{A}^N, \mathcal{B}^N)$ be a vector of $N \in \mathbb{N}$ input random variables fulfilling Assumption 2.1 and $\mathbf{X}_u := (X_{i_1}, \ldots, X_{i_{|u|}})^T : (\Omega^u, \mathcal{F}^u) \to (\mathbb{A}^u, \mathcal{B}^u)$, $\emptyset \neq u \subseteq \{1, \ldots, N\}$, be a subvector of $\mathbf{X}$. Then $\{P_{u,\mathbf{j}_u}(\mathbf{x}_u) : |\mathbf{j}_u| = l, \mathbf{j}_u \in \mathbb{N}^{|u|}\}$, the set of multivariate orthogonal polynomials of degree $l$, $|u| \leq l < \infty$, consistent with the probability measure $f_{\mathbf{X}_u}(\mathbf{x}_u)d\mathbf{x}_u$, is a basis of $\mathcal{Z}_l^u$.

*Proof.* Under Assumption 2.1, orthogonal polynomials consistent with the probability measure $f_{\mathbf{X}_u}(\mathbf{x}_u)d\mathbf{x}_u$ exist. Denote by $\mathbf{P}_{u,l} = (P_{u,l}^{(1)}, \ldots, P_{u,l}^{(M_{u,l})})^T$ a column vector of the elements of $\{P_{u,\mathbf{j}_u}(\mathbf{X}_u) : |\mathbf{j}_u| = l, \mathbf{j}_u \in \mathbb{N}^{|u|}\}$ arranged according to some monomial order. Let $\mathbf{a}_{u,l}^T = (a_{u,l}^{(1)}, \ldots, a_{u,l}^{(M_{u,l})})$ be a row vector comprising some constants $a_{u,l}^{(j)} \in \mathbb{R}$, $j = 1, \ldots, M_{u,l}$. Set $\mathbf{a}_{u,l}^T \mathbf{P}_{u,l} = 0$. Multiply both sides of the equality from the right by $\mathbf{P}_{u,l}^T$, integrate with respect to the measure $f_{\mathbf{X}_u}(\mathbf{x}_u)d\mathbf{x}_u$ over $\mathbb{A}^u$, and apply transposition to obtain

$$\mathbf{G}_{u,l}\mathbf{a}_{u,l} = \mathbf{0}, \tag{3.19}$$

where $\mathbf{G}_{u,l} = \mathbb{E}[\mathbf{P}_{u,l}\mathbf{P}_{u,l}^T]$ is an $M_{u,l} \times M_{u,l}$ matrix with its $(p,q)$th element

$$G_{u,l}^{(pq)} = \int_{\mathbb{A}^u} P_{u,l}^{(p)}(\mathbf{x}_u) P_{u,l}^{(q)}(\mathbf{x}_u) f_{\mathbf{X}_u}(\mathbf{x}_u) d\mathbf{x}_u = \mathbb{E}\left[P_{u,l}^{(p)}(\mathbf{X}_u) P_{u,l}^{(q)}(\mathbf{X}_u)\right]$$

representing the covariance between two elements of $\mathbf{P}_{u,l}$. According to Proposition 3.3, any two distinct polynomials from $\{P_{u,\mathbf{j}_u}(\mathbf{x}_u) : |\mathbf{j}_u| = l, \mathbf{j}_u \in \mathbb{N}^{|u|}\}$ are orthogonal, meaning that $\mathbb{E}[P_{u,l}^{(p)} P_{u,l}^{(q)}]$ is zero if $p \neq q$ and positive and finite if $p = q$. Consequently, $\mathbf{G}_{u,l}$ is a diagonal, positive-definite matrix and hence invertible. Therefore, (3.19) yields $\mathbf{a}_{u,l} = \mathbf{0}$, proving linear independence of the elements of $\mathbf{P}_{u,l}$ or the set $\{P_{u,\mathbf{j}_u}(\mathbf{x}_u) : |\mathbf{j}_u| = l, \mathbf{j}_u \in \mathbb{N}^{|u|}\}$. Furthermore, the dimension of $\mathcal{Z}_l^u$, which is $M_{u,l}$, matches exactly the number of elements of the aforementioned set. Therefore, the spanning set $\{P_{u,\mathbf{j}_u}(\mathbf{x}_u) : |\mathbf{j}_u| = l, \mathbf{j}_u \in \mathbb{N}^{|u|}\}$ forms a basis of $\mathcal{Z}_l^u$. ∎

*Proposition 3.6.* Let $\mathbf{X} := (X_1, \ldots, X_N)^T : (\Omega, \mathcal{F}) \to (\mathbb{A}^N, \mathcal{B}^N)$ be a vector of $N \in \mathbb{N}$ input random variables fulfilling both Assumptions 2.1 and 2.2 and $\mathbf{X}_u := (X_{i_1}, \ldots, X_{i_{|u|}})^T : (\Omega^u, \mathcal{F}^u) \to (\mathbb{A}^u, \mathcal{B}^u)$, $\emptyset \neq u \subseteq \{1, \ldots, N\}$, be a subvector of $\mathbf{X}$. Consistent with the probability measure $f_{\mathbf{X}_u}(\mathbf{x}_u)d\mathbf{x}_u$, let $\{P_{u,\mathbf{j}_u}(\mathbf{x}_u) : |\mathbf{j}_u| = l, \mathbf{j}_u \in \mathbb{N}^{|u|}\}$, the set of multivariate orthogonal polynomials of degree $l$, $|u| \leq l < \infty$, be a basis of $\mathcal{Z}_l^u$. Then the set of polynomials from the orthogonal sum

$$\mathbf{1} \oplus \bigoplus_{\emptyset \neq u \subseteq \{1,\ldots,N\}} \bigoplus_{l=|u|}^{\infty} \mathrm{span}\{P_{u,\mathbf{j}_u}(\mathbf{x}_u) : |\mathbf{j}_u| = l, \mathbf{j}_u \in \mathbb{N}^{|u|}\}$$

is dense in $L^2(\mathbb{A}^N, \mathcal{B}^N, f_{\mathbf{X}}d\mathbf{x})$. Moreover,

$$L^2(\mathbb{A}^N, \mathcal{B}^N, f_{\mathbf{X}}d\mathbf{x}) = \overline{\mathbf{1} \oplus \bigoplus_{\emptyset \neq u \subseteq \{1,\ldots,N\}} \bigoplus_{l=|u|}^{\infty} \mathcal{Z}_l^u} \tag{3.20}$$

*where the overline denotes set closure.*

*Proof.* Under Assumption 2.1, orthogonal polynomials exist. According to Theorem 3.4 of Ernst *et al.* [10], which exploits Assumption 2.2, the polynomial space $\Pi^{\{i\}} = \mathbb{R}[x_i]$ is dense in $L^2(\mathbb{A}^{\{i\}}, \mathcal{B}^{\{i\}}, f_{X_i} dx_i)$. Now use Theorem 4 of Petersen [19], which asserts that if, for $p \geq 1$ and all $i = 1, \ldots, N$, $\Pi^{\{i\}}$ is dense in $L^p(\mathbb{A}^{\{i\}}, \mathcal{B}^{\{i\}}, f_{X_i} dx_i)$, then so is $\Pi^N = \mathbb{R}[x_1, \ldots, x_N]$ in $L^p(\mathbb{A}^N, \mathcal{B}^N, f_{\mathbf{X}} d\mathbf{x})$. Therefore, the set of polynomials from the orthogonal sum, which is equal to $\Pi^N$ as per (3.13), is dense in $L^2(\mathbb{A}^N, \mathcal{B}^N, f_{\mathbf{X}} d\mathbf{x})$. Including the limit points of the orthogonal sum yields (3.20). ∎

**4. Polynomial Dimensional Decomposition.** Let $y(\mathbf{X}) := y(X_1, \ldots, X_N)$ be a real-valued, square-integrable output random variable defined on the same probability space $(\Omega, \mathcal{F}, \mathbb{P})$. The vector space $L^2(\Omega, \mathcal{F}, \mathbb{P})$ is a Hilbert space such that

$$\mathbb{E}\left[y^2(\mathbf{X})\right] := \int_\Omega y^2(\mathbf{X}(\omega)) d\mathbb{P}(\omega) = \int_{\mathbb{A}^N} y^2(\mathbf{x}) f_{\mathbf{X}}(\mathbf{x}) d\mathbf{x} < \infty$$

with inner product

$$(y(\mathbf{X}), z(\mathbf{X}))_{L^2(\Omega, \mathcal{F}, \mathbb{P})} := \int_\Omega y(\mathbf{X}(\omega)) z(\mathbf{X}(\omega)) d\mathbb{P}(\omega) = \int_{\mathbb{A}^N} y(\mathbf{x}) z(\mathbf{x}) f_{\mathbf{X}}(\mathbf{x}) d\mathbf{x} =: (y(\mathbf{x}), z(\mathbf{x}))_{f_{\mathbf{X}} d\mathbf{x}}$$

and norm

$$\|y(\mathbf{X})\|_{L^2(\Omega, \mathcal{F}, \mathbb{P})} := \sqrt{(y(\mathbf{X}), y(\mathbf{X}))_{L^2(\Omega, \mathcal{F}, \mathbb{P})}} = \sqrt{(y(\mathbf{x}), y(\mathbf{x}))_{f_{\mathbf{X}} d\mathbf{x}}} =: \|y(\mathbf{x})\|_{f_{\mathbf{X}} d\mathbf{x}}.$$

It is elementary to show that $y(\mathbf{X}(\omega)) \in L^2(\Omega, \mathcal{F}, \mathbb{P})$ if and only if $y(\mathbf{x}) \in L^2(\mathbb{A}^N, \mathcal{B}^N, f_{\mathbf{X}} d\mathbf{x})$.

**4.1. ADD.** The ADD, expressed by the recursive form [17, 22]

$$(4.1a) \qquad y(\mathbf{X}) = y_\emptyset + \sum_{\emptyset \neq u \subseteq \{1, \ldots, N\}} y_u(\mathbf{X}_u),$$

$$(4.1b) \qquad y_\emptyset = \int_{\mathbb{A}^N} y(\mathbf{x}) f_{\mathbf{X}}(\mathbf{x}) d\mathbf{x},$$

$$(4.1c) \qquad y_u(\mathbf{X}_u) = \int_{\mathbb{A}^{N-|u|}} y(\mathbf{X}_u, \mathbf{x}_{-u}) f_{\mathbf{X}_{-u}}(\mathbf{x}_{-u}) d\mathbf{x}_{-u} - \sum_{v \subset u} y_v(\mathbf{X}_v),$$

is a finite, hierarchical expansion of $y$ in terms of its input variables with increasing dimensions, where $u \subseteq \{1, \cdots, N\}$ is a subset with the complementary set $-u = \{1, \cdots, N\} \setminus u$ and $y_u$ is a $|u|$-variate component function describing a constant or an $|u|$-variate interaction of $\mathbf{X}_u = (X_{i_1}, \cdots, X_{i_{|u|}})$ on $y$ when $|u| = 0$ or $|u| > 0$. Here, $(\mathbf{X}_u, \mathbf{x}_{-u})$ denotes an $N$-dimensional vector whose $i$th component is $X_i$ if $i \in u$ and $x_i$ if $i \notin u$. The summation in (4.1a) comprises $2^N - 1$ terms with each term depending on a group of variables indexed by a particular subset of $\{1, \cdots, N\}$. When $u = \emptyset$, the sum in (4.1c) vanishes, resulting in the expression of the constant function $y_\emptyset$ in (4.1b). When $u = \{1, \cdots, N\}$, the integration in the last line of (4.1c) is on the empty set, reproducing (4.1a) and hence finding the last function $y_{\{1, \cdots, N\}}$. Indeed, all component functions of $y$ can be obtained by interpreting literally (4.1c). This decomposition, first presented by Hoeffding [16] in relation to his seminal work on $U$-statistics, has been studied by many other researchers described by Efron and Stein [9], the author [22], and references cited therein.

The ADD can also be generated by tensorizing a univariate function space decomposition into its constant subspace and remainder, producing [14]

$$L^2(\mathbb{A}^N, \mathcal{B}^N, f_{\mathbf{X}} d\mathbf{x}) = \mathbf{1} \oplus \bigoplus_{\emptyset \neq u \subseteq \{1,\ldots,N\}} \mathcal{W}_u, \tag{4.2}$$

where

$$\mathcal{W}_u := \{y_u \in L^2(\mathbb{A}^u, \mathcal{B}^u, f_{\mathbf{X}_u} d\mathbf{x}_u) : \mathbb{E}[y_u(\mathbf{X}_u) y_v(\mathbf{X}_v)] = 0 \text{ if } u \neq v, v \subseteq \{1,\ldots,N\}\}$$

is an ADD subspace comprising $|u|$-variate component functions of $y$. However, the subspaces $\mathcal{W}_u$, $\emptyset \neq u \subseteq \{1,\ldots,N\}$, are in general infinite-dimensional; therefore, further discretization of $\mathcal{W}_u$ is necessary. For instance, by introducing measure-consistent orthogonal polynomial basis discussed in Section 3, a component function $y_u(\mathbf{X}_u) \in \mathcal{W}_u$ can be expressed as a linear combination of these basis functions. Indeed, comparing (3.20) and (4.2) yields the closure of an orthogonal decomposition of

$$\mathcal{W}_u = \overline{\bigoplus_{l=|u|}^{\infty} \mathcal{Z}_l^u} \tag{4.3}$$

into polynomial spaces $\mathcal{Z}_l^u$, $|u| \leq l < \infty$. The result is a polynomial refinement of ADD, which is commonly referred to as PDD.

**4.2. PDD.** The PDD of a square-integrable random variable $y(\mathbf{X}) \in L^2(\Omega, \mathcal{F}, \mathbb{P})$ is simply the expansion of $y(\mathbf{X})$ with respect to a complete, hierarchically ordered, orthonormal polynomial basis of $L^2(\Omega, \mathcal{F}, \mathbb{P})$. There are at least two ways to explain PDD: a polynomial variant of ADD and a dimension-wise orthogonal polynomial expansion.

**4.2.1. Polynomial variant of ADD.** The first approach, explained by the author in a prior work [20], involves the following two steps: (1) expand the ANOVA component function

$$y_u(\mathbf{X}_u) \sim \sum_{\mathbf{j}_u \in \mathbb{N}^{|u|}} C_{u,\mathbf{j}_u} \Psi_{u,\mathbf{j}_u}(\mathbf{X}_u) \tag{4.4}$$

in terms of the basis of $\mathcal{W}_u$, which originally stems from the basis of $\mathcal{Z}_l^u$, $|u| \leq l < \infty$, with

$$C_{u,\mathbf{j}_u} = \int_{\mathbb{A}^{|u|}} y_u(\mathbf{X}_u) \Psi_{u,\mathbf{j}_u}(\mathbf{x}_u) f_{\mathbf{X}_u}(\mathbf{x}_u) d\mathbf{x}_u, \quad \emptyset \neq u \subseteq \{1,\ldots,N\}, \mathbf{j}_u \in \mathbb{N}^{|u|}, \tag{4.5}$$

representing the associated expansion coefficients; and (2) apply (4.1c) to (4.5) and exploit orthogonal properties of the basis. The end result is the PDD [20] of

$$y(\mathbf{X}) \sim y_\emptyset + \sum_{\emptyset \neq u \subseteq \{1,\ldots,N\}} \sum_{\mathbf{j}_u \in \mathbb{N}^{|u|}} C_{u,\mathbf{j}_u} \Psi_{u,\mathbf{j}_u}(\mathbf{X}_u), \tag{4.6}$$

where, eventually,

$$C_{u,\mathbf{j}_u} = \int_{\mathbb{A}^N} y(\mathbf{x}) \Psi_{u,\mathbf{j}_u}(\mathbf{x}_u) f_{\mathbf{X}}(\mathbf{x}) d\mathbf{x}. \tag{4.7}$$

Comparing (4.1) and (4.6), the connection between PDD and ADD is clearly palpable, where the former can be viewed as a polynomial variant of the latter. For instance, $C_{u,\mathbf{j}_u} \Psi_{u,\mathbf{j}_u}(\mathbf{X}_u)$ in (4.6) represents a $|u|$-variate, $|\mathbf{j}_u|$th-order PDD component function of $y(\mathbf{X})$, describing the $|\mathbf{j}_u|$th-order polynomial approximation of $y_u(\mathbf{X}_u)$. In addition, PDD inherits all desirable properties of ADD [20].

### 4.2.2. Dimension-wise orthogonal polynomial expansion.
The second approach entails polynomial expansion associated with the dimension-wise orthogonal splitting of polynomial spaces, as explained in Section 3.3. The latter approach has not been published elsewhere and is, therefore, formally presented here as Theorem 4.1.

Theorem 4.1. *Let $\mathbf{X} := (X_1, \ldots, X_N)^T : (\Omega, \mathcal{F}) \to (\mathbb{A}^N, \mathcal{B}^N)$ be a vector of $N \in \mathbb{N}$ input random variables fulfilling Assumptions 2.1 and 2.2. For $\emptyset \neq u \subseteq \{1, \ldots, N\}$ and $\mathbf{X}_u := (X_{i_1}, \ldots, X_{i_{|u|}})^T : (\Omega^u, \mathcal{F}^u) \to (\mathbb{A}^u, \mathcal{B}^u)$, denote by $\{\Psi_{u,\mathbf{j}_u}(\mathbf{X}_u) : \mathbf{j}_u \in \mathbb{N}^{|u|}\}$ the set of multivariate orthonormal polynomials consistent with the probability measure $f_{\mathbf{X}_u}(\mathbf{x}_u) d\mathbf{x}_u$. Then*

(1) *any random variable $y(\mathbf{X}) \in L^2(\Omega, \mathcal{F}, \mathbb{P})$ can be hierarchically expanded as a Fourier-like series, referred to as the PDD of*

$$
\begin{aligned}
y(\mathbf{X}) \;\sim\; & y_\emptyset + \sum_{\emptyset \neq u \subseteq \{1,\ldots,N\}} \sum_{l=|u|}^{\infty} \sum_{\substack{\mathbf{j}_u \in \mathbb{N}^{|u|} \\ |\mathbf{j}_u| = l}} C_{u,\mathbf{j}_u} \Psi_{u,\mathbf{j}_u}(\mathbf{X}_u) \\
= \; & y_\emptyset + \sum_{\emptyset \neq u \subseteq \{1,\ldots,N\}} \sum_{\mathbf{j}_u \in \mathbb{N}^{|u|}} C_{u,\mathbf{j}_u} \Psi_{u,\mathbf{j}_u}(\mathbf{X}_u),
\end{aligned}
\tag{4.8}
$$

*where the expansion coefficients $y_\emptyset \in \mathbb{R}$ and $C_{u,\mathbf{j}_u} \in \mathbb{R}$, $\emptyset \neq u \subseteq \{1, \ldots, N\}$, $\mathbf{j}_u \in \mathbb{N}^{|u|}$, are defined by*

$$
y_\emptyset := \mathbb{E}\left[y(\mathbf{X})\right] := \int_{\mathbb{A}^N} y(\mathbf{x}) f_{\mathbf{X}}(\mathbf{x}) d\mathbf{x},
\tag{4.9}
$$

$$
C_{u,\mathbf{j}_u} := \mathbb{E}\left[y(\mathbf{X}) \Psi_{u,\mathbf{j}_u}(\mathbf{X}_u)\right] := \int_{\mathbb{A}^N} y(\mathbf{x}) \Psi_{u,\mathbf{j}_u}(\mathbf{x}_u) f_{\mathbf{X}}(\mathbf{x}) d\mathbf{x};
\tag{4.10}
$$

*and*

(2) *the PDD of $y(\mathbf{X}) \in L^2(\Omega, \mathcal{F}, \mathbb{P})$ converges to $y(\mathbf{X})$ in mean-square; furthermore, the PDD converges in probability and in distribution.*

*Proof.* Under Assumptions 2.1 and 2.2, a complete infinite set of multivariate orthogonal polynomials in $\mathbf{x}_u$ consistent with the probability measure $f_{\mathbf{X}_u}(\mathbf{x}_u) d\mathbf{x}_u$ exists. From Proposition 3.6 and the fact that orthonormality is merely scaling, the set of polynomials from the orthogonal sum

$$
\mathbf{1} \oplus \bigoplus_{\emptyset \neq u \subseteq \{1,\ldots,N\}} \bigoplus_{l=|u|}^{\infty} \operatorname{span}\{\Psi_{u,\mathbf{j}_u}(\mathbf{x}_u) : |\mathbf{j}_u| = l, \mathbf{j}_u \in \mathbb{N}^{|u|}\} = \Pi^N
\tag{4.11}
$$

is also dense in $L^2(\mathbb{A}^N, \mathcal{B}^N, f_{\mathbf{X}} d\mathbf{x})$. Therefore, any square-integrable random variable $y(\mathbf{X})$ can be expanded as shown in (4.8). Combining the two inner sums of the expansion forms the equality in the second line of (4.8).

From the denseness, one has the Bessel's inequality [7]

$$
\mathbb{E}\left[y_\emptyset + \sum_{\emptyset \neq u \subseteq \{1,\ldots,N\}} \sum_{\mathbf{j}_u \in \mathbb{N}^{|u|}} C_{u,\mathbf{j}_u} \Psi_{u,\mathbf{j}_u}(\mathbf{X}_u)\right]^2 \leq \mathbb{E}\left[y^2(\mathbf{X})\right],
\tag{4.12}
$$

proving that the PDD converges in mean-square or $L^2$. To determine the limit of convergence, invoke again Proposition 3.6, which implies that the set on the left side of (4.11) is complete in $L^2(\mathbb{A}^N, \mathcal{B}^N, f_\mathbf{X} d\mathbf{x})$. Therefore, Bessel's inequality becomes an equality

$$(4.13) \quad \mathbb{E}\left[y_\emptyset + \sum_{\emptyset \neq u \subseteq \{1,\ldots,N\}} \sum_{\mathbf{j}_u \in \mathbb{N}^{|u|}} C_{u,\mathbf{j}_u} \Psi_{u,\mathbf{j}_u}(\mathbf{X}_u)\right]^2 = \mathbb{E}\left[y^2(\mathbf{X})\right],$$

known as the Parseval identity [7] for a multivariate orthonormal system, for every random variable $y(\mathbf{X}) \in L^2(\Omega, \mathcal{F}, \mathbb{P})$. Furthermore, as the PDD converges in mean-square, it does so in probability. Moreover, as the expansion converges in probability, it also converges in distribution.

Finally, to find the expansion coefficients, define a second moment

$$(4.14) \quad e_{\mathrm{PDD}} := \mathbb{E}\left[y(\mathbf{X}) - y_\emptyset - \sum_{\emptyset \neq v \subseteq \{1,\ldots,N\}} \sum_{\mathbf{k}_v \in \mathbb{N}^{|v|}} C_{v,\mathbf{k}_v} \Psi_{v,\mathbf{k}_v}(\mathbf{X}_v)\right]^2$$

of the difference between $y(\mathbf{X})$ and its full PDD. Differentiate both sides of (4.14) with respect to $y_\emptyset$ and $C_{u,\mathbf{j}_u}$, $\emptyset \neq u \subseteq \{1,\ldots,N\}$, $\mathbf{j}_u \in \mathbb{N}^{|u|}$, to write

$$(4.15) \quad \begin{aligned} \frac{\partial e_{\mathrm{PDD}}}{\partial y_\emptyset} &= \frac{\partial}{\partial y_\emptyset} \mathbb{E}\left[y(\mathbf{X}) - y_\emptyset - \sum_{\emptyset \neq v \subseteq \{1,\ldots,N\}} \sum_{\mathbf{k}_v \in \mathbb{N}^{|v|}} C_{v,\mathbf{k}_v} \Psi_{v,\mathbf{k}_v}(\mathbf{X}_v)\right]^2 \\ &= \mathbb{E}\left[\frac{\partial}{\partial y_\emptyset}\left\{y(\mathbf{X}) - y_\emptyset - \sum_{\emptyset \neq v \subseteq \{1,\ldots,N\}} \sum_{\mathbf{k}_v \in \mathbb{N}^{|v|}} C_{v,\mathbf{k}_v} \Psi_{v,\mathbf{k}_v}(\mathbf{X}_v)\right\}^2\right] \\ &= 2\mathbb{E}\left[\left\{y_\emptyset + \sum_{\emptyset \neq v \subseteq \{1,\ldots,N\}} \sum_{\mathbf{k}_v \in \mathbb{N}^{|v|}} C_{v,\mathbf{k}_v} \Psi_{v,\mathbf{k}_v}(\mathbf{X}_v) - y(\mathbf{X})\right\} \times 1\right] \\ &= 2\left\{y_\emptyset - \mathbb{E}\left[y(\mathbf{X})\right]\right\} \end{aligned}$$

and

$$(4.16) \quad \begin{aligned} \frac{\partial e_{\mathrm{PDD}}}{\partial C_{u,\mathbf{j}_u}} &= \frac{\partial}{\partial C_{u,\mathbf{j}_u}} \mathbb{E}\left[y(\mathbf{X}) - y_\emptyset - \sum_{\emptyset \neq v \subseteq \{1,\ldots,N\}} \sum_{\mathbf{k}_v \in \mathbb{N}^{|v|}} C_{v,\mathbf{k}_v} \Psi_{v,\mathbf{k}_v}(\mathbf{X}_v)\right]^2 \\ &= \mathbb{E}\left[\frac{\partial}{\partial C_{u,\mathbf{j}_u}}\left\{y(\mathbf{X}) - y_\emptyset - \sum_{\emptyset \neq v \subseteq \{1,\ldots,N\}} \sum_{\mathbf{k}_v \in \mathbb{N}^{|v|}} C_{v,\mathbf{k}_v} \Psi_{v,\mathbf{k}_v}(\mathbf{X}_v)\right\}^2\right] \\ &= 2\mathbb{E}\left[\left\{y_\emptyset + \sum_{\emptyset \neq v \subseteq \{1,\ldots,N\}} \sum_{\mathbf{k}_v \in \mathbb{N}^{|v|}} C_{v,\mathbf{k}_v} \Psi_{v,\mathbf{k}_v}(\mathbf{X}_v) - y(\mathbf{X})\right\} \Psi_{u,\mathbf{j}_u}(\mathbf{X}_u)\right] \\ &= 2\left\{C_{u,\mathbf{j}_u} - \mathbb{E}\left[y(\mathbf{X}) \Psi_{u,\mathbf{j}_u}(\mathbf{X}_u)\right]\right\}. \end{aligned}$$

Here, the second, third, and last lines of both (4.15) and (4.16) are obtained by interchanging the differential and expectation operators, performing the differentiation, swapping the expectation and summation operators and applying Corollary 3.4, respectively. The interchanges are permissible as the infinite sum is convergent as demonstrated in the preceding paragraph. Setting $\partial e_{\mathrm{PDD}}/\partial y_\emptyset = 0$ in (4.15) and $\partial e_{\mathrm{PDD}}/\partial C_{u,\mathbf{j}_u} = 0$ in (4.16) yields (4.9) and (4.10), respectively, completing the proof. ∎

The expressions of the expansion coefficients can also be derived by simply replacing $y(\mathbf{X})$ in (4.9) and (4.10) with the full PDD and then using Corollary 3.4. In contrast, the proof given here demonstrates that the PDD coefficients are determined optimally.

It should be emphasized that the function $y$ must be square-integrable for the mean-square and other convergences to hold. However, the rate of convergence depends on the smoothness of the function. The smoother the function, the faster the convergence. If the function is a polynomial, then its PDD exactly reproduces the function. These results can be easily proved using classical approximation theory.

A related expansion, known by the name of RS-HDMR [18], also involves orthogonal polynomials in connection with ADD. However, the existence, convergence, and approximation quality of the expansion, including its behavior for infinitely many input variables, have not been reported.

**4.3. Truncation.** The full PDD contains an infinite number of orthonormal polynomials or coefficients. In practice, the number must be finite, meaning that PDD must be truncated. However, there are multiple ways to perform the truncation. A straightforward approach adopted in this work entails (1) keeping all polynomials in at most $0 \leq S \leq N$ variables, thereby retaining the degrees of interaction among input variables less than or equal to $S$ and (2) preserving polynomial expansion orders (total) less than or equal to $S \leq m < \infty$. The result is an $S$-variate, $m$th-order PDD approximation[2]

$$
\begin{aligned}
y_{S,m}(\mathbf{X}) &= y_\emptyset + \sum_{s=1}^{S} \sum_{l=s}^{m} \sum_{\substack{\emptyset \neq u \subseteq \{1,\ldots,N\} \\ |u|=s}} \sum_{\substack{\mathbf{j}_u \in \mathbb{N}^{|u|} \\ |\mathbf{j}_u|=l}} C_{u,\mathbf{j}_u} \Psi_{u,\mathbf{j}_u}(\mathbf{X}_u) \\
&= y_\emptyset + \sum_{\substack{\emptyset \neq u \subseteq \{1,\ldots,N\} \\ 1 \leq |u| \leq S}} \sum_{\substack{\mathbf{j}_u \in \mathbb{N}^{|u|} \\ |u| \leq |\mathbf{j}_u| \leq m}} C_{u,\mathbf{j}_u} \Psi_{u,\mathbf{j}_u}(\mathbf{X}_u)
\end{aligned}
\tag{4.17}
$$

of $y(\mathbf{X})$, containing

$$
L_{S,m} = 1 + \sum_{s=1}^{S} \binom{N}{s}\binom{m}{s}
\tag{4.18}
$$

number of expansion coefficients including $y_\emptyset$. It is important to clarify a few things about the truncated PDD proposed. First, a different truncation with respect to the polynomial expansion order based on $\infty$-norm as opposed to 1-norm, that is, $\|\mathbf{j}_u\|_\infty \leq m$, was employed in prior works [20, 21, 24]. Therefore, comparing (4.17) and (4.18) with the existing truncation, if it is desired, should be done with care. Having said this, the proposed truncation has one advantage over the existing one: a direct comparison with a truncated PCE is possible; this will be further explained in the forthcoming sections. Second, the right side of (4.17) contains sums of at most $S$-dimensional orthonormal polynomials, representing at most $S$-variate PDD component functions of $y$. Therefore, the term "$S$-variate" used for the PDD approximation should be interpreted in the context of including at most $S$-degree interaction of input variables, even though $y_{S,m}$ is strictly an $N$-variate function. Third, when $S = 0$, $y_{0,m} = y_\emptyset$ for any $m$ as the outer sums of (4.17) vanish. Finally, when $S \to N$ and $m \to \infty$, $y_{S,m}$ converges to $y$ in the mean-square sense, generating a hierarchical and convergent sequence of PDD approximations. Readers interested in an adaptive version of PDD, where the truncation parameters are automatically chosen, are directed to the work of Yadav and Rahman [30], including an application to design optimization [25].

---

[2]The nouns *degree* and *order* associated with PDD or orthogonal polynomials are used synonymously in the paper.

It is natural to ask about the approximation quality of (4.17). Since the set of polynomials from the orthogonal sum in (4.11) is complete in $L^2(\mathbb{A}^N, \mathcal{B}^N, f_\mathbf{X} d\mathbf{x})$, the truncation error $y(\mathbf{X}) - y_{S,m}(\mathbf{X})$ is orthogonal to any element of the subspace from which $y_{S,m}(\mathbf{X})$ is chosen, as demonstrated below.

**Proposition 4.2.** *For any $y(\mathbf{X}) \in L^2(\Omega, \mathcal{F}, \mathbb{P})$, let $y_{S,m}(\mathbf{X})$ be its $S$-variate, $m$th-order PDD approximation. Then the truncation error $y(\mathbf{X}) - y_{S,m}(\mathbf{X})$ is orthogonal to the subspace*

$$\Pi^N_{S,m} := \mathbf{1} \oplus \bigoplus_{\substack{\emptyset \neq u \subseteq \{1,\ldots,N\} \\ 1 \leq |u| \leq S}} \bigoplus_{\substack{\mathbf{j}_u \in \mathbb{N}^{|u|} \\ |u| \leq |\mathbf{j}_u| \leq m}} span\{\Psi_{u,\mathbf{j}_u}(\mathbf{X}_u) : \mathbf{j}_u \in \mathbb{N}^{|u|}\} \subseteq L^2(\Omega, \mathcal{F}, \mathbb{P}), \tag{4.19}$$

*comprising all polynomials in $\mathbf{X}$ with the degree of interaction at most $S$ and order at most $m$, including constants. Moreover, $\mathbb{E}[y(\mathbf{X}) - y_{S,m}(\mathbf{X})]^2 \to 0$ as $S \to N$ and $m \to \infty$.*

*Proof.* Let

$$\bar{y}_{S,m}(\mathbf{X}) := \bar{y}_\emptyset + \sum_{\substack{\emptyset \neq v \subseteq \{1,\ldots,N\} \\ 1 \leq |v| \leq S}} \sum_{\substack{\mathbf{k}_v \in \mathbb{N}^{|v|} \\ |v| \leq |\mathbf{k}_v| \leq m}} \bar{C}_{v,\mathbf{k}_v} \Psi_{v,\mathbf{k}_v}(\mathbf{X}_v), \tag{4.20}$$

with arbitrary expansion coefficients $\bar{y}_\emptyset$ and $\bar{C}_{v,\mathbf{k}_v}$, be any element of the subspace $\Pi^N_{S,m}$ of $L^2(\Omega, \mathcal{F}, \mathbb{P})$ described by (4.19). Then
(4.21)
$$\begin{aligned}
&\mathbb{E}\left[\{y(\mathbf{X}) - y_{S,m}(\mathbf{X})\} \bar{y}_{S,m}(\mathbf{X})\right] \\
&= \mathbb{E}\left[\left\{\sum_{\substack{\emptyset \neq u \subseteq \{1,\ldots,N\} \\ 1 \leq |u| \leq S}} \sum_{\substack{\mathbf{j}_u \in \mathbb{N}^{|u|} \\ m+1 \leq |\mathbf{j}_u| < \infty}} C_{u,\mathbf{j}_u} \Psi_{u,\mathbf{j}_u}(\mathbf{X}_u) + \sum_{\substack{\emptyset \neq u \subseteq \{1,\ldots,N\} \\ S+1 \leq |u| \leq N}} \sum_{\substack{\mathbf{j}_u \in \mathbb{N}^{|u|} \\ |u| \leq |\mathbf{j}_u| < \infty}} C_{u,\mathbf{j}_u} \Psi_{u,\mathbf{j}_u}(\mathbf{X}_u)\right\}\right. \\
&\quad \times \left.\left\{\bar{y}_\emptyset + \sum_{\substack{\emptyset \neq v \subseteq \{1,\ldots,N\} \\ 1 \leq |v| \leq S}} \sum_{\substack{\mathbf{k}_v \in \mathbb{N}^{|v|} \\ |v| \leq |\mathbf{k}_v| \leq m}} \bar{C}_{v,\mathbf{k}_v} \Psi_{v,\mathbf{k}_v}(\mathbf{X}_v)\right\}\right] \\
&= 0,
\end{aligned}$$

where the last line follows from Corollary 3.4, proving the first part of the proposition. For the latter part, the Pythagoras theorem yields

$$\mathbb{E}[\{y(\mathbf{X}) - y_{S,m}(\mathbf{X})\}^2] + \mathbb{E}[y^2_{S,m}(\mathbf{X})] = \mathbb{E}[y^2(\mathbf{X})]. \tag{4.22}$$

From Theorem 4.1, $\mathbb{E}[y^2_{S,m}(\mathbf{X})] \to \mathbb{E}[y^2(\mathbf{X})]$ as $S \to N$ and $m \to \infty$. Therefore, $\mathbb{E}[\{y(\mathbf{X}) - y_{S,m}(\mathbf{X})\}^2] \to 0$ as $S \to N$ and $m \to \infty$. ∎

The second part of Proposition 4.2 entails $L^2$ convergence, which is the same as the mean-square convergence described in Theorem 4.1. However, an alternative route is chosen for the proof of the proposition. Besides, Proposition 4.2 implies that the PDD approximation is optimal as it recovers the best approximation from the subspace $\Pi^N_{S,m}$, as described by Corollary 4.3.

Corollary 4.3. *Let* $\Pi_{S,m}^N$ *in* (4.19) *define the subspace of all polynomials in* $\mathbf{X}$ *with the degree of interaction at most* $S$ *and order at most* $m$, *including constants. Then the* $S$-*variate,* $m$th*-order PDD approximation* $y_{S,m}(\mathbf{X})$ *of* $y(\mathbf{X}) \in L^2(\Omega, \mathcal{F}, \mathbb{P})$ *is the best approximation in the sense that*

$$\text{(4.23)} \qquad \mathbb{E}\left[y(\mathbf{X}) - y_{S,m}(\mathbf{X})\right]^2 = \inf_{\bar{y}_{S,m} \in \Pi_{S,m}^N} \mathbb{E}\left[y(\mathbf{X}) - \bar{y}_{S,m}(\mathbf{X})\right]^2.$$

*Proof.* Consider two elements $y_{S,m}(\mathbf{X})$ and $\bar{y}_{S,m}(\mathbf{X})$ of the subspace $\Pi_{S,m}^N$, where the former is the $S$-variate, $m$th-order PDD approximation of $y(\mathbf{X})$ with the expansion coefficients defined by (4.9) and (4.10) and the latter is any $S$-variate, $m$th-order polynomial function, described by (4.20), with arbitrary chosen expansion coefficients. From Proposition 4.2, the truncation error $y(\mathbf{X}) - y_{S,m}(\mathbf{X})$ is orthogonal to both $y_{S,m}(\mathbf{X})$ and $\bar{y}_{S,m}(\mathbf{X})$ and is, therefore, orthogonal to their linear combinations, yielding

$$\mathbb{E}\left[\{y(\mathbf{X}) - y_{S,m}(\mathbf{X})\}\{y_{S,m}(\mathbf{X}) - \bar{y}_{S,m}(\mathbf{X})\}\right] = 0.$$

Consequently,

$$\text{(4.24)} \qquad \begin{aligned} \mathbb{E}\left[y(\mathbf{X}) - \bar{y}_{S,m}(\mathbf{X})\right]^2 &= \mathbb{E}\left[y(\mathbf{X}) - y_{S,m}(\mathbf{X})\right]^2 + \mathbb{E}\left[y_{S,m}(\mathbf{X}) - \bar{y}_{S,m}(\mathbf{X})\right]^2 \\ &\geq \mathbb{E}\left[y(\mathbf{X}) - y_{S,m}(\mathbf{X})\right]^2, \end{aligned}$$

as the second expectation on the right side of the first line of (4.24) is non-negative, thereby proving the mean-square optimality of the $S$-variate, $m$th-order PDD approximation. ∎

The motivations behind ADD- and PDD-derived approximations are the following. In a practical setting, the function $y(\mathbf{X})$, fortunately, has an effective dimension [3] much lower than $N$, meaning that the right side of (4.1a) can be effectively approximated by a sum of lower-dimensional component functions $y_u$, $|u| \ll N$, but still maintaining all random variables $\mathbf{X}$ of a high-dimensional uncertainty quantification problem. Furthermore, an $S$-variate, $m$th-order PDD approximation is grounded on a fundamental conjecture known to be true in many real-world uncertainty quantification problems: given a high-dimensional function $y$, its $|u|$-variate, $|\mathbf{j}_u|$th-order PDD component function $C_{u,\mathbf{j}_u}\Psi_{u,\mathbf{j}_u}(\mathbf{X}_u)$, where $S+1 \leq |u| \leq N$ and $m+1 \leq |\mathbf{j}_u| < \infty$, is small and hence negligible, leading to an accurate low-variate, low-order approximation of $y$. The computational complexity of a truncated PDD is polynomial, as opposed to exponential, thereby alleviating the curse of dimensionality to a substantial extent. Although PCE contains the same orthogonal polynomials, a recent work on random eigenvalue analysis of dynamic systems reveals markedly higher convergence rate of the PDD approximation than the PCE approximation [24].

**4.4. Output statistics and other probabilistic characteristics.** The $S$-variate, $m$th-order PDD approximation $y_{S,m}(\mathbf{X})$ can be viewed as a surrogate of $y(\mathbf{X})$. Therefore, relevant probabilistic characteristics of $y(\mathbf{X})$, including its first two moments and probability density function, if it exists, can be estimated from the statistical properties of $y_{S,m}(\mathbf{X})$.

Applying the expectation operator on $y_{S,m}(\mathbf{X})$ and $y(\mathbf{X})$ in (4.17) and (4.8) and imposing Corollary 3.4, their means

$$\text{(4.25)} \qquad \mathbb{E}\left[y_{S,m}(\mathbf{X})\right] = \mathbb{E}\left[y(\mathbf{X})\right] = y_\emptyset$$

are the same and independent of $S$ and $m$. Therefore, the PDD truncated for any values of $0 \le S \le N$ and $S \le m < \infty$ yields the exact mean. Nonetheless, $\mathbb{E}[y_{S,m}(\mathbf{X})]$ will be referred to as the $S$-variate, $m$th-order PDD approximation of the mean of $y(\mathbf{X})$.

Applying the expectation operator again, this time on $[y_{S,m}(\mathbf{X}) - y_\emptyset]^2$ and $[y(\mathbf{X}) - y_\emptyset]^2$, and employing Corollary 3.4 results in the variances

$$\text{(4.26)} \qquad \text{var}\,[y_{S,m}(\mathbf{X})] = \sum_{\substack{\emptyset \ne u \subseteq \{1,\ldots,N\} \\ 1 \le |u| \le S}} \sum_{\substack{\mathbf{j}_u \in \mathbb{N}^{|u|} \\ |u| \le |\mathbf{j}_u| \le m}} C^2_{u,\mathbf{j}_u}$$

and

$$\text{(4.27)} \qquad \text{var}\,[y(\mathbf{X})] = \sum_{\emptyset \ne u \subseteq \{1,\ldots,N\}} \sum_{\mathbf{j}_u \in \mathbb{N}^{|u|}} C^2_{u,\mathbf{j}_u}$$

of $y_{S,m}(\mathbf{X})$ and $y(\mathbf{X})$, respectively. Again, $\text{var}[y_{S,m}(\mathbf{X})]$ will be referred to as the $S$-variate, $m$th-order PDD approximation of the variance of $y(\mathbf{X})$. Clearly, $\text{var}[y_{S,m}(\mathbf{X})]$ approaches $\text{var}[y(\mathbf{X})]$, the exact variance of $y(\mathbf{X})$, as $S \to N$ and $m \to \infty$.

Being convergent in probability and distribution, the probability density function of $y(\mathbf{X})$, if it exists, can also be estimated by that of $y_{S,m}(\mathbf{X})$. However, no analytical formula exists for the density function. In that case, the density can be estimated by sampling methods, such as Monte Carlo simulation (MCS) of $y_{S,m}(\mathbf{X})$. Such simulation should not be confused with crude MCS of $y(\mathbf{X})$, commonly used for producing benchmark results whenever possible. The crude MCS can be expensive or even prohibitive, particularly when the sample size needs to be very large for estimating tail probabilistic characteristics. In contrast, the MCS embedded in the PDD approximation requires evaluations of simple polynomial functions that describe $y_{S,m}$. Therefore, a relatively large sample size can be accommodated in the PDD approximation even when $y$ is expensive to evaluate.

**4.5. Infinitely many input variables.** In many fields, such as uncertainty quantification, information theory, and stochastic process, functions depending on a countable sequence $\{X_i\}_{i \in \mathbb{N}}$ of input random variables need to be considered [15]. Under certain assumptions, PDD is still applicable as in the case of finitely many random variables, as demonstrated by the following proposition.

Proposition 4.4. *Let $\{X_i\}_{i \in \mathbb{N}}$ be a countable sequence of input random variables defined on the probability space $(\Omega, \mathcal{F}_\infty, \mathbb{P})$, where $\mathcal{F}_\infty := \sigma(\{X_i\}_{i \in \mathbb{N}})$ is the associated $\sigma$-algebra generated. If the sequence $\{X_i\}_{i \in \mathbb{N}}$ satisfies Assumptions 2.1 and 2.2, then the PDD of $y(\{X_i\}_{i \in \mathbb{N}}) \in L^2(\Omega, \mathcal{F}_\infty, \mathbb{P})$, where $y : \mathbb{A}^\mathbb{N} \to \mathbb{R}$, converges to $y(\{X_i\}_{i \in \mathbb{N}})$ in mean-square. Moreover, the PDD converges in probability and in distribution.*

*Proof.* According to Proposition 3.6, $\Pi^N$ is dense in $L^2(\mathbb{A}^N, \mathcal{B}^N, f_\mathbf{X} d\mathbf{x})$ and hence in $L^2(\Omega, \mathcal{F}_N, \mathbb{P})$ for every $N \in \mathbb{N}$, where $\mathcal{F}_N := \sigma(\{X_i\}_{i=1}^N)$ is the associated $\sigma$-algebra generated by $\{X_i\}_{i=1}^N$. Here, with a certain abuse of notation, $\Pi^N$ is used as a set of polynomial functions of both real variables $\mathbf{x}$ and random variables $\mathbf{X}$. Now, apply Theorem 3.8 of Ernst et al. [10], which says that if $\Pi^N$ is dense in $L^2(\Omega, \mathcal{F}_N, \mathbb{P})$ for every $N \in \mathbb{N}$, then

$$\Pi^\infty := \bigcup_{N=1}^\infty \Pi^N,$$

a subspace of $L^2(\Omega, \mathcal{F}_\infty, \mathbb{P})$, is also dense in $L^2(\Omega, \mathcal{F}_\infty, \mathbb{P})$. But, using (4.11),

$$\begin{aligned}\Pi^\infty &= \bigcup_{N=1}^\infty \mathbf{1} \oplus \bigoplus_{\emptyset \neq u \subseteq \{1,\ldots,N\}} \bigoplus_{l=|u|}^\infty \text{span}\{\Psi_{u,\mathbf{j}_u} : |\mathbf{j}_u| = l, \mathbf{j}_u \in \mathbb{N}^{|u|}\} \\ &= \mathbf{1} \oplus \bigoplus_{\emptyset \neq u \subseteq \mathbb{N}} \bigoplus_{l=|u|}^\infty \text{span}\{\Psi_{u,\mathbf{j}_u} : |\mathbf{j}_u| = l, \mathbf{j}_u \in \mathbb{N}^{|u|}\},\end{aligned}$$

demonstrating that the set of polynomials from the orthogonal sum in the last line is dense in $L^2(\Omega, \mathcal{F}_\infty, \mathbb{P})$. Therefore, the PDD of $y(\{X_i\}_{i\in\mathbb{N}}) \in L^2(\Omega, \mathcal{F}_\infty, \mathbb{P})$ converges to $y(\{X_i\}_{i\in\mathbb{N}})$ in mean-square. Since the mean-square convergence is stronger than the convergence in probability or in distribution, the latter modes of convergence follow readily. ∎

**5. Polynomial Chaos Expansion.** In contrast to the dimension-wise splitting of polynomial spaces in PDD, a degree-wise orthogonal splitting of polynomial spaces results in PCE. The latter decomposition is briefly summarized here as PCE will be compared with PDD in the next section.

**5.1. Degree-wise orthogonal decomposition of polynomial spaces.** Let $\mathbf{j} := \mathbf{j}_{\{1,\ldots,N\}} = (j_1, \ldots, j_N) \in \mathbb{N}_0^N$, $j_i \in \mathbb{N}_0$, $i \in \{1,\ldots,N\}$, define an $N$-dimensional multi-index. For $\mathbf{x} = (x_1, \ldots, x_N) \in \mathbb{A}^N \subseteq \mathbb{R}^N$, a monomial in the variables $x_1, \ldots, x_N$ is the product $\mathbf{x}^\mathbf{j} = x_1^{j_1} \cdots x_N^{j_N}$ and has a total degree $|\mathbf{j}| = j_1 + \cdots + j_N$. Denote by

$$\Pi_p^N := \text{span}\{\mathbf{x}^\mathbf{j} : 0 \leq |\mathbf{j}| \leq p, \mathbf{j} \in \mathbb{N}_0^N\}, \ p \in \mathbb{N}_0,$$

the space of real polynomials in $\mathbf{x}$ of degree at most $p$. Let $\mathcal{V}_0^N := \Pi_0^N = \text{span}\{1\}$ be the space of constant functions. For each $1 \leq l < \infty$, denote by $\mathcal{V}_l^N \subset \Pi_l^N$ the space of orthogonal polynomials of degree exactly $l$ that are orthogonal to all polynomials in $\Pi_{l-1}^N$, that is,

$$\mathcal{V}_l^N := \{P \in \Pi_l^N : (P,Q)_{f_\mathbf{x} d\mathbf{x}} = 0 \ \forall Q \in \Pi_{l-1}^N\}, \ 1 \leq l < \infty.$$

From Section 3, with $u = \{1, \ldots, N\}$ in mind, select $\{P_\mathbf{j}(\mathbf{x}) : |\mathbf{j}| = l, \mathbf{j} \in \mathbb{N}_0^N\} \subset \mathcal{V}_l^N$ to be a basis of $\mathcal{V}_l^N$. Each basis function $P_\mathbf{j}(\mathbf{x})$ is a multivariate orthogonal polynomial in $\mathbf{x}$ of degree $|\mathbf{j}|$. Obviously,

$$\mathcal{V}_l^N = \text{span}\{P_\mathbf{j}(\mathbf{x}) : |\mathbf{j}| = l, \mathbf{j} \in \mathbb{N}_0^N\}, \ 0 \leq l < \infty.$$

According to (3.7) with $u = \{1, \ldots, N\}$, $P_\mathbf{j}(\mathbf{x})$ is orthogonal to $P_\mathbf{k}(\mathbf{x})$ whenever $\mathbf{j} \neq \mathbf{k}$. Therefore, any two polynomial subspaces $\mathcal{V}_l^N$ and $\mathcal{V}_r^N$, where $0 \leq l, r < \infty$, are orthogonal whenever $l \neq r$. In consequence, there exists another orthogonal decomposition of

$$(5.1) \quad \Pi^N = \bigoplus_{l \in \mathbb{N}_0} \mathcal{V}_l^N = \bigoplus_{l \in \mathbb{N}_0} \text{span}\{P_\mathbf{j}(\mathbf{x}) : |\mathbf{j}| = l, \mathbf{j} \in \mathbb{N}_0^N\} = \text{span}\{P_\mathbf{j}(\mathbf{x}) : \mathbf{j} \in \mathbb{N}_0^N\}.$$

Compared with (3.13), (5.1) represents a degree-wise orthogonal decomposition of $\Pi^N$.

**5.2. PCE.** Given the degree-wise orthogonal decomposition of $\Pi^N$, the PCE of any square-integrable output random variable $y(\mathbf{X})$ is expressed by [5, 10, 28, 29]

$$(5.2) \quad y(\mathbf{X}) \sim \sum_{l=0}^\infty \sum_{\substack{\mathbf{j} \in \mathbb{N}_0^N \\ |\mathbf{j}|=l}} C_\mathbf{j} \Psi_\mathbf{j}(\mathbf{X}) = \sum_{\mathbf{j} \in \mathbb{N}_0^N} C_\mathbf{j} \Psi_\mathbf{j}(\mathbf{X}),$$

where $\{\Psi_{\mathbf{j}}(\mathbf{X}) : \mathbf{j} \in \mathbb{N}_0^N\}$ is an infinite set of measure-consistent multivariate orthonormal polynomials in $\mathbf{X}$ that can be obtained by scaling $P_{\mathbf{j}}$ in (3.14) and $C_{\mathbf{j}} \in \mathbb{R}$, $\mathbf{j} \in \mathbb{N}_0^N$, are the PCE expansion coefficients. Like PDD, the PCE of $y(\mathbf{X}) \in L^2(\Omega, \mathcal{F}, \mathbb{P})$ under Assumptions 2.1 and 2.2 also converges to $y(\mathbf{X})$ in mean-square, in probability, and in distribution.

Since the PCE of $y(\mathbf{X})$ in (5.2) is an infinite series, it must also be truncated in applications. A commonly adopted truncation is based on retaining orders of polynomials less than or equal to a specified total degree. In this regard, given $0 \leq p < \infty$, the $p$th-order PCE approximation of $y(\mathbf{X}) \in L^2(\Omega, \mathcal{F}, \mathbb{P})$ reads

$$(5.3) \quad y_p(\mathbf{X}) = \sum_{l=0}^{p} \sum_{\substack{\mathbf{j} \in \mathbb{N}_0^N \\ |\mathbf{j}|=l}} C_{\mathbf{j}} \Psi_{\mathbf{j}}(\mathbf{X}) = \sum_{\substack{\mathbf{j} \in \mathbb{N}_0^N \\ 0 \leq |\mathbf{j}| \leq p}} C_{\mathbf{j}} \Psi_{\mathbf{j}}(\mathbf{X}).$$

This kind of truncation is related to the total degree index set

$$\left\{ \mathbf{j} \in \mathbb{N}_0^N : \sum_{i=1}^{N} j_i \leq p \right\}$$

for defining the recovered multivariate polynomial space of a $p$th-order PCE approximation. Other kinds of truncation entail

$$\left\{ \mathbf{j} \in \mathbb{N}_0^N : \max_{i=1,\ldots,N} j_i \leq p \right\} \text{ and } \left\{ \mathbf{j} \in \mathbb{N}_0^N : \prod_{i=1}^{N} (j_i + 1) \leq p + 1 \right\},$$

describing the tensor product and hyperbolic cross index sets, respectively, to name just two. The total degree and tensor product index sets are common choices, although the latter one suffers from the curse of dimensionality, making it impractical for high-dimensional problems. The hyperbolic cross index set, originally introduced for approximating periodic functions by trigonometric polynomials [2], is relatively a new idea and has yet to receive widespread attention. All of these choices and possibly others, including their anisotropic versions, can be used for truncating PCE. In this work, however, only the total degree index set is used for the PCE approximation. This is consistent with the 1-norm of $\mathbf{j}_u$ used for truncating PDD in (4.17).

## 6. Error Analysis.

### 6.1. PDD error.
Define a second-moment error,

$$(6.1) \quad e_{S,m} := \mathbb{E}\left[y(\mathbf{X}) - y_{S,m}(\mathbf{X})\right]^2,$$

stemming from the $S$-variate, $m$th-order PDD approximation presented in the preceding section. Replacing $y(\mathbf{X})$ and $y_{S,m}(\mathbf{X})$ in (6.1) with the right sides of (4.8) and (4.17), respectively, produces

$$(6.2) \quad e_{S,m} = \sum_{s=1}^{S} \sum_{l=m+1}^{\infty} \sum_{\substack{\emptyset \neq u \subseteq \{1,\ldots,N\} \\ |u|=s}} \sum_{\substack{\mathbf{j}_u \in \mathbb{N}^{|u|} \\ |\mathbf{j}_u|=l}} C_{u,\mathbf{j}_u}^2 + \sum_{s=S+1}^{N} \sum_{l=s}^{\infty} \sum_{\substack{\emptyset \neq u \subseteq \{1,\ldots,N\} \\ |u|=s}} \sum_{\substack{\mathbf{j}_u \in \mathbb{N}^{|u|} \\ |\mathbf{j}_u|=l}} C_{u,\mathbf{j}_u}^2,$$

where the second term vanishes expectedly when $S = N$ as the lower limit of the outer sum exceeds the upper limit. In (6.2), the first term of the PDD error is due to the truncation of polynomial expansion orders involving interactive effects of at most $S$ variables, whereas the second term of the PDD error is contributed by ignoring the interactive effects of larger than $S$ variables. Obviously, the error for a general function $y$ depends on which expansion coefficients decay and how they decay with respect to $S$ and $m$. Nonetheless, the error decays monotonically with respect to $S$ and/or $m$ as stated in Proposition 6.1. Other than that, nothing more can be said about the PDD error.

**Proposition 6.1.** *For a general function $y$, $e_{S+i,m+j} \leq e_{S,m}$, where $0 \leq S < N$, $S \leq m < \infty$, and $i$ and $j$ are equal to either 0 or 1, but not both equal to 0.*

*Proof.* Setting $i = 1$, $j = 0$ and using (6.2),

$$e_{S+1,m} - e_{S,m} = \sum_{l=m+1}^{\infty} \sum_{\substack{\emptyset \neq u \subseteq \{1,\ldots,N\} \\ |u|=S+1}} \sum_{\substack{\mathbf{j}_u \in \mathbb{N}^{|u|} \\ |\mathbf{j}_u|=l}} C_{u,\mathbf{j}_u}^2 - \sum_{l=S+1}^{\infty} \sum_{\substack{\emptyset \neq u \subseteq \{1,\ldots,N\} \\ |u|=S+1}} \sum_{\substack{\mathbf{j}_u \in \mathbb{N}^{|u|} \\ |\mathbf{j}_u|=l}} C_{u,\mathbf{j}_u}^2 \leq 0,$$

where the inequality to *zero* in the last line results from the fact that, as $S \leq m$, the first term is smaller than or equal to the second term. Similarly, setting $i = 0$, $j = 1$ and using (6.2),

$$e_{S,m+1} - e_{S,m} = -\sum_{s=1}^{S} \sum_{\substack{\emptyset \neq u \subseteq \{1,\ldots,N\} \\ |u|=s}} \sum_{\substack{\mathbf{j}_u \in \mathbb{N}^{|u|} \\ |\mathbf{j}_u|=m+1}} C_{u,\mathbf{j}_u}^2 \leq 0.$$

Finally, setting $i = 1$, $j = 1$,

$$e_{S+1,m+1} - e_{S,m} = e_{S+1,m+1} - e_{S,m+1} + e_{S,m+1} - e_{S,m} \leq 0,$$

as $e_{S+1,m+1} - e_{S,m+1} \leq 0$ and $e_{S,m+1} - e_{S,m} \leq 0$. ∎

**Corollary 6.2.** *For a general function $y$, $e_{S',m'} \leq e_{S,m}$ whenever $S' \geq S$ and $m' \geq m$.*

In practice, the effects of interaction among input variables and polynomial expansion order become increasingly weaker as $|u|$ and $|\mathbf{j}_u|$ grow. In this case, $C_{u,\mathbf{j}_u}^2$, which is equal to the variance of $C_{u,\mathbf{j}_u} \Psi_{u,\mathbf{j}_u}(\mathbf{X}_u)$, decreases with $|u|$ and $|\mathbf{j}_u|$. Given the rates at which $C_{u,\mathbf{j}_u}^2$ decreases with $|u|$ and $|\mathbf{j}_u|$, a question arises on how fast does $e_{S,m}$ decay with respect to $S$ and $m$. Proposition 6.3, Corollary 6.4, and subsequent discussions provide a few insights.

**Proposition 6.3.** *For a class of functions $y$, assume that $C_{u,\mathbf{j}_u}^2$, $\emptyset \neq u \subseteq \{1,\cdots,N\}$, $\mathbf{j}_u \in \mathbb{N}^{|u|}$, attenuates according to $C_{u,\mathbf{j}_u}^2 \leq C p_1^{-|u|} p_2^{-|\mathbf{j}_u|}$, where $C > 0$, $p_1 > 1$, and $p_2 > 1$ are three real-valued constants. Then it holds that*

$$\text{var}[y(\mathbf{X})] \leq C \left[ \left\{ \frac{1 + p_1(p_2 - 1)}{p_1(p_2 - 1)} \right\}^N - 1 \right] \tag{6.3}$$

*and*

$$e_{S,m} \leq C \left[ \sum_{s=1}^{S} \sum_{l=m+1}^{\infty} \binom{N}{s} \binom{l-1}{s-1} p_1^{-s} p_2^{-l} + \sum_{s=S+1}^{N} \sum_{l=s}^{\infty} \binom{N}{s} \binom{l-1}{s-1} p_1^{-s} p_2^{-l} \right]. \tag{6.4}$$

*Proof.* With the recognition that

$$\#\{\emptyset \neq u \subseteq \{1,\cdots,N\} : |u| = s\} = \binom{N}{s}, \quad \#\{\mathbf{j}_u \in \mathbb{N}^{|u|} : |\mathbf{j}_u| = l\} = \binom{l-1}{|u|-1},$$

use $C_{u,\mathbf{j}_u}^2 \leq Cp_1^{-|u|}p_2^{-|\mathbf{j}_u|}$ in (4.27) and (6.2) to obtain (6.3) and (6.4). ∎

COROLLARY 6.4. *For the function class described in Proposition 6.3, $e_{S+i,m+j} < e_{S,m}$, where $0 \leq S < N$, $S \leq m < \infty$, and $i$ and $j$ are equal to either 0 or 1, but not both equal to 0.*

According to Corollary 6.4, $e_{S,m}$ decays strictly monotonically with respect to $S$ and/or $m$ for any rate parameters $p_1 > 1$ and $p_2 > 1$. When the equality holds in (6.3) and (6.4) from Proposition 6.3, Figure 1, comprising three subfigures, presents three sets of plots of the relative error, $e_{S,m}/\text{var}[y(\mathbf{X})]$, against $m$ for five distinct values of $S = 1, 2, 3, 5, 9$. These subfigures, each obtained for $N = 20$, correspond to three distinct cases of the values of $p_1$ and $p_2$: (1) $p_1 = 500$, $p_2 = 5$; (2) $p_1 = 5$, $p_2 = 500$; and (3) $p_1 = 500$, $p_2 = 500$. In all cases, the error for a given $S$ decays first with respect to $m$, and then levels off at a respective limit when $m$ is sufficiently large. The limits get progressively smaller when $S$ increases as expected. However, the magnitude of this behavior depends on the rates at which the expansion coefficient attenuates with respect to the degree of interaction and polynomial expansion order. When $p_1 > p_2$, as in case 1 [Figure 1 (top)], the error for a given $S$ decays slowly with respect to $m$ due to a relatively weaker attenuation rate associated with the polynomial expansion order. The trend reverses when the attenuation rate becomes stronger and reaches the condition $p_1 < p_2$, as in case 2 [Figure 1 (middle)]. For larger values of $S$, for example, $S = 5$ or $9$, the respective limits are significantly lower in case 2 than in case 1. When the attenuation rates are the same and large, as in case 3 [Figure 1 (bottom)], the decay rate of error accelerates substantially.

**6.2. Relationship between PDD and PCE.** Since PDD and PCE share the same orthonormal polynomials, they are related. Indeed, the relationship was first studied by Rahman and Yadav [24], who determined that any one of the two infinite series from PDD and PCE defined by (4.8) and (5.2) can be rearranged to derive the other. In other words, the PDD can also be viewed as a reshuffled PCE and vice versa. However, due to a strong connection to ADD endowed with a desired hierarchical structure, PDD merits its own appellation. More importantly, the PDD and PCE when truncated are not the same. In fact, two important observations stand out prominently. First, the terms in the PCE approximation are organized with respect to the order of polynomials. In contrast, the PDD approximation is structured with respect to the degree of interaction between a finite number of random variables. Therefore, significant differences may exist regarding the accuracy, efficiency, and convergence properties of their truncated sum. Second, if a stochastic response is highly nonlinear, but contains rapidly diminishing interactive effects of multiple random variables, the PDD approximation is expected to be more effective than the PCE approximation. This is because the lower-variate terms of the PDD approximation can be just as nonlinear by selecting appropriate values of $m$ in (4.17). In contrast, many more terms and expansion coefficients are required to be included in the PCE approximation to capture such high nonlinearity. In this work, a theoretical comparison between PDD and PCE in the context of error analysis, not studied in prior works, is presented.

For error analysis, it is convenient to write a PCE approximation in terms of a PDD approximation. Indeed, there exists a striking result connecting PCE with PDD approximations, as explained in Proposition 6.3.

PROPOSITION 6.5. *Let $y_p(\mathbf{X})$ and $y_{S,m}(\mathbf{X})$ be the pth-order PCE approximation and S-variate, mth-order PDD approximation of $y(\mathbf{X}) \in L^2(\Omega, \mathcal{F}, \mathbb{P})$, respectively, where $0 \leq S \leq N$, $S \leq m < \infty$, and $0 \leq p < \infty$. Then the pth-order PCE approximation and the*

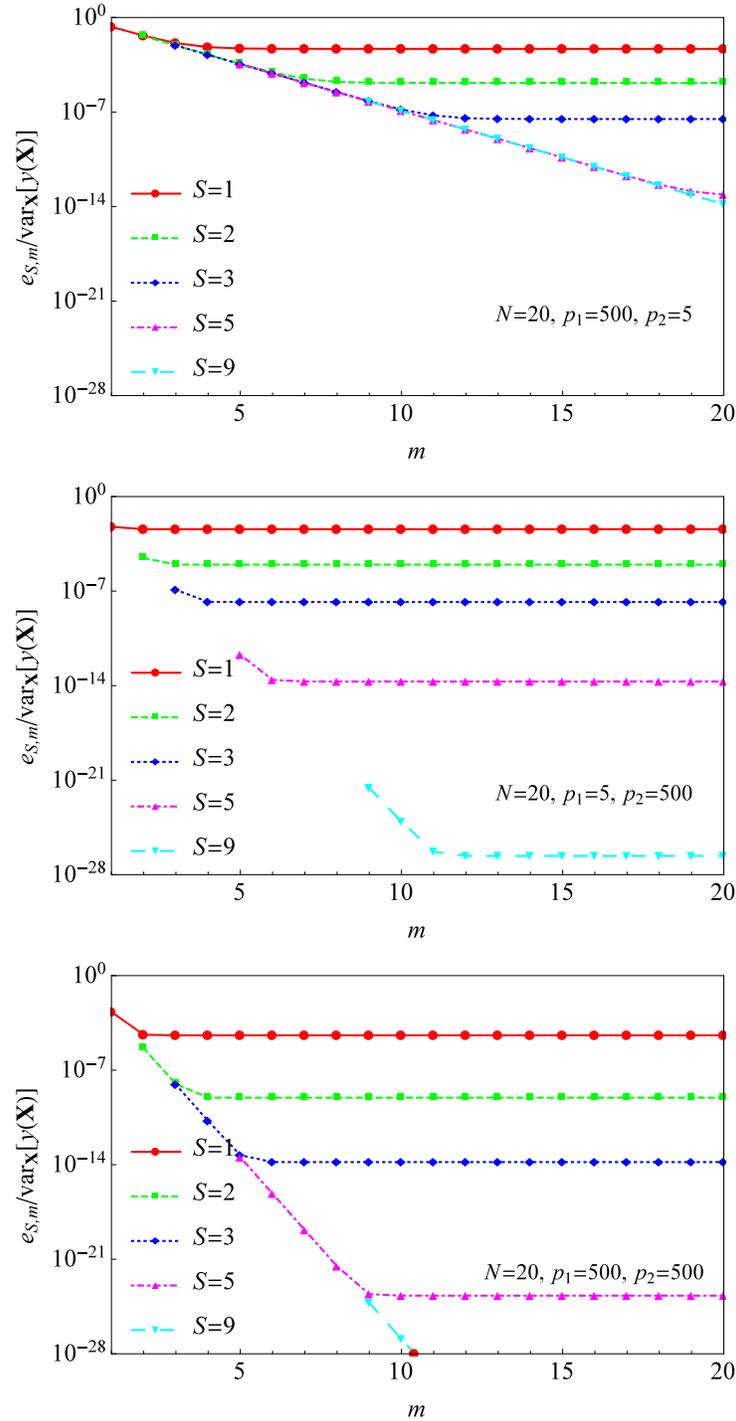

**Figure 1.** *PDD errors for various attenuation rates of the expansion coefficients; (top) $p_1 = 500$, $p_2 = 5$; (middle) $p_1 = 5$, $p_2 = 500$; (bottom) $p_1 = 500$, $p_2 = 500$.*

$(p \wedge N)$-variate, $p$th-order PDD approximation are the same, that is,

$$(6.5) \qquad y_p(\mathbf{X}) = y_{p \wedge N, p}(\mathbf{X}),$$

where $y_{0,0}(\mathbf{X}) = y_\emptyset$ and $p \wedge N$ denotes the minimum of $p$ and $N$.

*Proof.* According to Rahman and Yadav [24], the right side of (5.3) can be reshuffled, resulting in a long form of the PCE approximation, expressed by

$$(6.6) \quad y_p(\mathbf{X}) = y_\emptyset + \sum_{s=1}^{N} \left[ \underbrace{\sum_{i_1=1}^{N-s+1} \cdots \sum_{i_s=i_{s-1}+1}^{N}}_{s \text{ sums}} \underbrace{\sum_{j_{i_1}=1}^{p-s+1} \cdots \sum_{j_{i_s}=1}^{p-s+1}}_{s \text{ sums}; j_{i_1}+\cdots+j_{i_s} \leq p} C_{\{i_1 \cdots i_s\},(j_{i_1} \cdots j_{i_s})} \prod_{q=1}^{s} \psi_{i_q j_{i_q}}(X_{i_q}) \right],$$

in terms of the PDD expansion coefficients. Note that, depending on the condition $p \leq N$ or $p \geq N$, at most $p$-dimensional or $N$-dimensional sums survive in (6.6), meaning that the $p$th-order PCE approximation retains effects of at most $(p \wedge N)$-degree interaction and at most $p$th-order polynomial expansion order. Accordingly, the compact form of the PCE approximation can be written as

$$(6.7) \qquad y_p(\mathbf{X}) = y_\emptyset + \sum_{s=1}^{p \wedge N} \sum_{l=s}^{p} \sum_{\substack{\emptyset \neq u \subseteq \{1,\ldots,N\} \\ |u|=s}} \sum_{\substack{\mathbf{j}_u \in \mathbb{N}^{|u|} \\ |\mathbf{j}_u|=l}} C_{u,\mathbf{j}_u} \Psi_{u,\mathbf{j}_u}(\mathbf{X}_u) = y_{p \wedge N, p}(\mathbf{X}),$$

completing the proof. ∎

Using Proposition 6.5, the number of expansion coefficients, say, $L_p$ associated with the $p$th-order PCE approximation can be calculated from that required by the $(p \wedge N)$-variate, $p$th-order PDD approximation. Accordingly, setting $S = p \wedge N$ and $m = p$ in (4.18),

$$(6.8) \qquad L_p = L_{p \wedge N, p} = 1 + \sum_{s=1}^{p \wedge N} \binom{N}{s} \binom{p}{s} = \frac{(N+p)!}{N! p!}$$

with the last expression commonly found in the PCE literature [29]. The advantage of (6.7) over (5.3) is obvious: the PDD coefficients, once determined, can be reused for the PCE approximation and subsequent error analysis, thereby sidestepping calculations of the PCE coefficients.

### 6.3. PDD vs. PCE errors. 
Define another second-moment error,

$$(6.9) \qquad e_p := \mathbb{E}\left[y(\mathbf{X}) - y_p(\mathbf{X})\right]^2,$$

resulting from the $p$th-order PCE approximation $y_p(\mathbf{X})$ of $y(\mathbf{X})$. Using Proposition 6.5, $e_p = e_{p \wedge N, p}$, meaning that the PCE error analysis can be conducted using the PDD approximation.

**Proposition 6.6.** *For a general function $y$, let $e_{S,m}$ and $e_p$ denote the PDD and PCE errors defined by (6.1) and (6.9), respectively. Given a truncation parameter $0 \leq p < \infty$ of the PCE approximation, if the truncation parameters of the PDD approximation are chosen such that $p \wedge N \leq S \leq N$ and $p \vee S \leq m < \infty$, then*

$$(6.10) \qquad e_{S,m} \leq e_p,$$

where $p \vee S$ denotes the maximum of $p$ and $S$.

*Proof.* The result follows from Propositions 6.1 and 6.5, and Corollary 6.2. ∎

Proposition 6.6 aids in selecting appropriate truncation parameters to contrast the second-moment errors due to PDD and PCE approximations. However, the proposition does not say anything about the computational effort. Proposition 6.7 and subsequent discussion explain the relationship between computational effort and error committed by both PDD and PCE approximations for a special class of functions.

**Proposition 6.7.** *For a special class of functions $y$, assume that $C^2_{u,\mathbf{j}_u}$, $\emptyset \neq u \subseteq \{1, \cdots, N\}$, $\mathbf{j}_u \in \mathbb{N}^{|u|}$, diminishes according to $C^2_{u,\mathbf{j}_u} \leq C p_1^{-|u|} p_2^{-|\mathbf{j}_u|}$, where $C > 0$, $p_1 > 1$, and $p_2 > 1$ are three real-valued constants. Then it holds that*

$$(6.11) \quad e_p \leq C \left[ \sum_{s=1}^{p \wedge N} \sum_{l=p+1}^{\infty} \binom{N}{s}\binom{l-1}{s-1} p_1^{-s} p_2^{-l} + \sum_{s=p \wedge N+1}^{N} \sum_{l=s}^{\infty} \binom{N}{s}\binom{l-1}{s-1} p_1^{-s} p_2^{-l} \right].$$

*Proof.* Replacing $S$ and $m$ in (6.4) with $p \wedge N$ and $p$, respectively, obtains the result. ∎

Theoretically, the numbers of expansion coefficients required by the PDD and PCE approximations can be used to compare their respective computational efforts. Table 1 presents for $N = 20$ the requisite numbers of expansion coefficients when PDD is truncated at $S = 1, 2, 3, 5, 9$ and $m = 1 - 20$, and when PCE is truncated at $p = 1 - 20$. They are calculated using (4.18) and (6.8) for PDD and PCE approximations, respectively. According to Table 1, the growth of the number of expansion coefficients in PCE is steeper than that in PDD. The growth rate increases markedly when the polynomial expansion order is large. This is primarily because a PCE approximation is solely dictated by a single truncation parameter $p$, which controls the largest polynomial expansion order preserved, but not the degree of interaction independently. In contrast, two different truncation parameters $S$ and $m$ are involved in a PDD approximation, affording a greater flexibility in retaining the largest degree of interaction and largest polynomial expansion order. In consequence, the numbers of expansion coefficients and hence the computational efforts by the PDD and PCE approximations can vary appreciably.

**Table 1**
*Growth of expansion coefficients in the PDD and PCE approximations*

| | $L_{S,m}$ | | | | | |
|---|---|---|---|---|---|---|
| $m$ or $p$ | $S=1$ | $S=2$ | $S=3$ | $S=5$ | $S=9$ | $L_p$ |
| 1 | 21 | | | | | 21 |
| 2 | 41 | 231 | | | | 231 |
| 3 | 61 | 631 | 1771 | | | 1771 |
| 5 | 101 | 2001 | 13,401 | 53,130 | | 53,130 |
| 9 | 181 | 7021 | 102,781 | 2,666,755 | 10,015,005 | 10,015,005 |
| 12 | 241 | 12,781 | 263,581 | 14,941,024 | 211,457,454 | 225,792,840 |
| 15 | 301 | 20,251 | 538,951 | 53,710,888 | 2,397,802,638 | 3,247,943,160 |
| 20 | 401 | 36,501 | 1,336,101 | 265,184,142 | 51,855,874,642 | 137,846,528,820 |

Using the equalities in (6.3), (6.4), and (6.11), Figure 2 depicts how the relative PDD error, $e_{S,m}/\text{var}[y(\mathbf{X})]$, and the relative PCE error $e_p/\text{var}[y(\mathbf{X})]$, vary with respect to the number of expansion coefficients required for $N = 20$. Again, the three preceding cases of the attenuation rates $p_1$ and $p_2$ with respect to the degree of interaction and polynomial expansion order are studied. In all cases, the PDD or PCE errors decay with respect to $S$, $m$, and $p$ as expected. However, in the PDD approximation, the error for a fixed $S$ may decline even further by increasing $m$, whereas no such possibility exists in the PCE approximation. This behavior is pronounced in case 1, that is, when $p_1 > p_2$ [Figure 2 (top)]. For example, in case 1, the bivariate, sixth-order PDD approximation ($S = 2$, $m = 6$) achieves a relative error of $8.54 \times 10^{-5}$ employing only 2971 expansion coefficients. In contrast, to match the same-order error, the sixth-order PCE approximation ($p = 6$) is needed, committing a relative error of $7.15 \times 10^{-5}$ at the cost of 230,230 expansion coefficients. Therefore, the PDD approximation is substantially more economical than the PCE approximation for a similar accuracy. However, when $p_1 > p_2$, as in case 2 [Figure 2 (middle)], the computational advantage of PDD over PCE approximations disappears as the attenuation rate associated with the polynomial expansion order is dominant over that associated with the degree of interaction. Nonetheless, in case 2, the $S$-variate, $m$th-order PDD approximation with the lowest $m$ possible cannot commit more error than the $m$th-order PCE approximation for the same computational effort. Finally, when the attenuation rates are the same, as in case 3 [Figure 2 (bottom)], the PDD approximation is still more computationally efficient than the PCE approximation. For instance, the trivariate, fifth-order PDD ($S = 3$, $m = 5$) and fifth-order PCE ($p = 6$) approximations require 13,401 and 53,130 expansion coefficients to commit the same-order errors of $5.07 \times 10^{-14}$ and $3.51 \times 10^{-14}$, respectively. But, unlike in case 1, an unnecessarily large polynomial expansion order may render the PDD approximation more expensive than required.

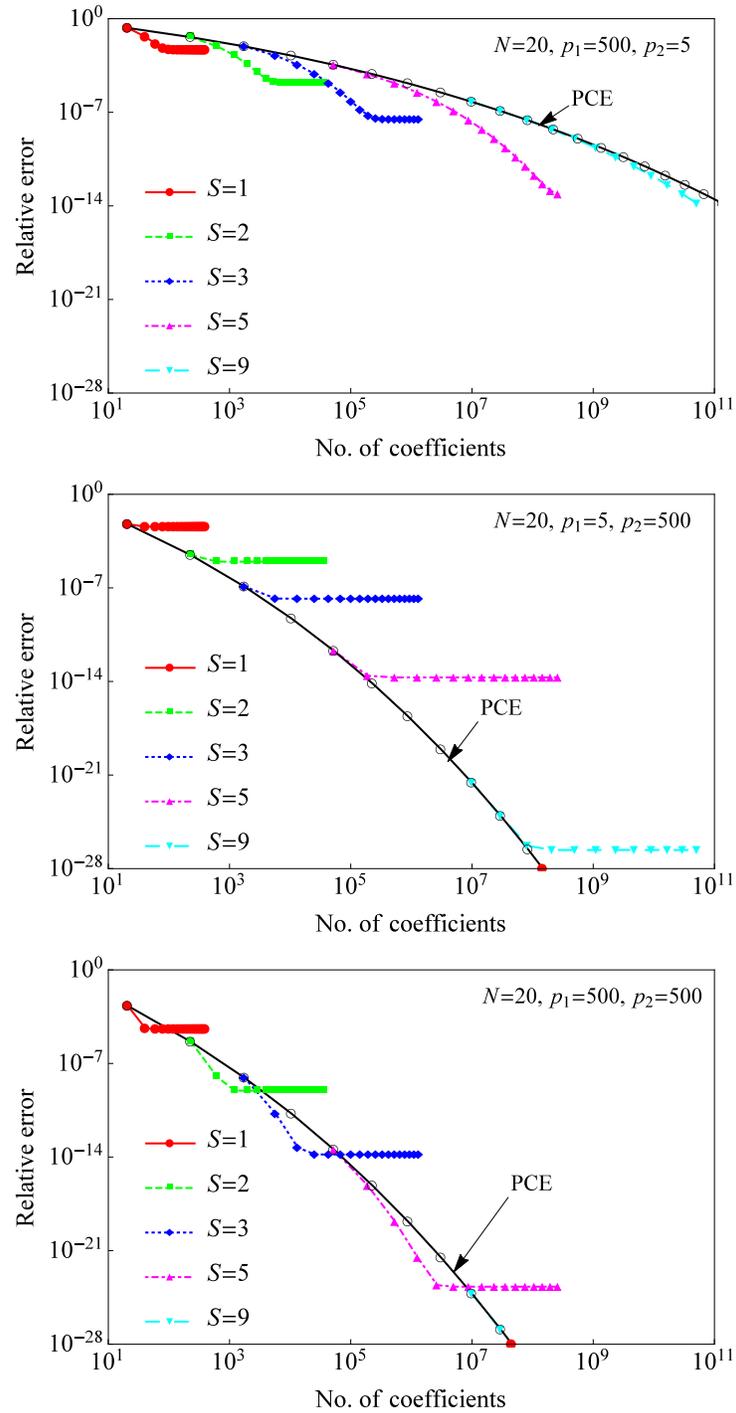

**Figure 2.** *PDD vs. PCE errors for various attenuation rates of the expansion coefficients; (top) $p_1 = 500$, $p_2 = 5$; (middle) $p_1 = 5$, $p_2 = 500$; (bottom) $p_1 = 500$, $p_2 = 500$.*

Readers should take note that the comparative error analyses reported here are limited to PDD and PCE approximations derived from truncations according to the total degree index set. For other index sets, such as the tensor product and hyperbolic cross index sets, it would be intriguing to find whether a similar conclusion arises.

**7. Conclusion.** The fundamental mathematical properties of PDD, representing Fourier-like series expansion in terms of random orthogonal polynomials with increasing dimensions, were studied. A dimension-wise splitting of appropriate polynomial spaces into orthogonal subspaces, each spanned by measure-consistent orthogonal polynomials, was constructed, resulting in a polynomial refinement of ADD and eventually PDD. Under prescribed assumptions, the set of measure-consistent orthogonal polynomials was proved to form a complete basis of each subspace, leading to an orthogonal sum of such sets of basis functions, including the constant subspace, to span the space of all polynomials. In addition, the orthogonal sum is dense in a Hilbert space of square-integrable functions, leading to mean-square convergence of PDD to the correct limit, including for the case of infinitely many random variables. The optimality of PDD and the approximation quality due to truncation were demonstrated or discussed. From the second-moment error analysis of a general function of $1 \leq N < \infty$ random variables, given $0 \leq p < \infty$, the $(p \wedge N)$-variate, $p$th-order PDD approximation and $p$th-order PCE approximation are the same. Therefore, an $S$-variate, $m$th-order PDD approximation cannot commit a larger error than a $p$th-order PCE approximation if $p \wedge N \leq S \leq N$ and $p \vee S \leq m < \infty$. From the comparison of computational efforts, required to estimate with the same accuracy the variance of an output function entailing exponentially attenuating expansion coefficients, the PDD approximation can be substantially more economical than the PCE approximation.